\documentclass[11pt]{article}
\usepackage{amssymb,amsmath,latexsym}

\newtheorem{thm}{Theorem}[section]
\newtheorem{cor}[thm]{Corollary}
\newtheorem{lemma}[thm]{Lemma}
\newtheorem{prop}[thm]{Proposition}

\newtheorem{remark}[thm]{Remark}
\numberwithin{equation}{section}

\def\nn{\nonumber}

\def\pf{{\medskip\noindent {\bf Proof. }}}
\def\qed{{\hfill $\Box$ \bigskip}}

 \def\sE {{\cal E}} \def\sF {{\cal F}}
  
  \def\sL {{\cal L}}

  \def\bR {{\mathbb R}}

\def\R {{\mathbb R}}
\def\wt{\widetilde}

\def\E{{\mathbb E}}
\def\P{{\mathbb P}}

\def\bea{\begin{align*}}
\def\eea{\end{align*}}
\def\bee{\begin{equation}}
\def\eee{\end{equation}}

\def\eps{\varepsilon}

\makeatletter
\@addtoreset{equation}{section}

\makeatother

\begin{document}
\allowdisplaybreaks
\bibliographystyle{plain}

\title{\Large \bf
Dirichlet Heat Kernel Estimates for Subordinate Brownian Motions with Gaussian Components}

\author{{\bf Zhen-Qing Chen}\thanks{Research partially supported
by   NSF Grant   DMS-1206276  and NNSFC Grant 11128101.},
\quad {\bf Panki Kim}\thanks{This research was supported by Basic Science Research Program through the National Research Foundation of Korea(NRF) funded by the Ministry of Education, Science and Technology (0409-20120034).} \quad and \quad {\bf Renming
Song}\thanks{Research supported in part by a grant from the Simons
Foundation (208236).} }
 \date{(March 26, 2013)}

\maketitle

\begin{abstract}
In this paper, we derive
explicit sharp two-sided estimates
 for the Dirichlet heat kernels, in $C^{1,1}$ open sets
$D$ in  $\R^d$, of a large class of subordinate Brownian motions
with Gaussian components.
When $D$ is bounded, our sharp two-sided Dirichlet heat kernel estimates hold for all $t>0$.
Integrating the heat kernel estimates with respect to the time variable $t$, we obtain
sharp two-sided estimates for the Green functions, in bounded $C^{1,1}$ open sets, of
such subordinate Brownian motions with Gaussian components.
\end{abstract}

\bigskip
\noindent {\bf AMS 2000 Mathematics Subject Classification}: Primary
60J35, 47G20, 60J75; Secondary 47D07

\bigskip\noindent
{\bf Keywords and phrases}:
Subordinate Brownian motion, Dirichlet heat kernel, transition density,
Green function, exit time, L\'evy system,
boundary Harnack principle

\bigskip

\section{Introduction}

It is well-known that, for a second order elliptic differential
operator $\sL$ on $\bR^d$ satisfying some natural conditions, there
is a diffusion process $X$ on $\bR^d$ with $\sL$ as its
infinitesimal generator. The fundamental solution $p(t, x, y)$ of
$\partial_tu=\sL u$ (also called the heat kernel of $\sL$) is the
transition density of $X$.
Such relationship is also true for a large
class of Markov processes with discontinuous sample paths, which
constitute an important family of stochastic processes in
probability theory that have been widely used in various applications.
Thus obtaining sharp two-sided
estimates for $p(t, x, y)$ is a fundamental problem in both analysis
and probability theory.

Two-sided heat kernel estimates for diffusions on $\bR^d$ have a
long history and many beautiful results have been established. See
\cite{D1, Deb} and the references therein. But, due to the
complication near the boundary, two-sided estimates for the
transition density
(equivalently, the Dirichlet heat kernels) of killed
Brownian motion in a
connected open set $D$
 have been established
only recently. See \cite{D2, Deb, DS} for upper bound estimates and
\cite{Zq3} for the lower bound estimates of the Dirichlet heat
kernels in bounded $C^{1,1}$
connected open sets.
For discontinuous processes (or, non-local operators), the study of
their global heat kernel estimates started quite recently. See
\cite{CKK2, CKK3, CK, CK2, CK3} and the references therein.  See also \cite{C}
for a recent survey on this topic. The study of sharp
two-sided Dirichlet heat kernel estimates
for discontinuous processes is even more recent.
In \cite{CKS},
we obtained sharp two-sided estimates for the Dirichlet heat
kernel of the fractional Laplacian $\Delta^{\alpha /2}$ in any
$C^{1,1}$ open set $D$ with zero
exterior condition on $D^c$ (or equivalently, the transition density
function of the killed $\alpha$-stable process in $D$).

In the last few years, the
approach developed in \cite{CKS} has served as
a road map
for establishing sharp two-sided Dirichlet heat kernel estimates
for other purely discontinuous processes in open subsets of $\R^d$.
In \cite{CKS1, CKS6, CKS2}, the ideas of \cite{CKS}
were adapted and further developed to establish sharp two-sided
Dirichlet heat kernel
estimates of censored stable-like processes, mixed stable processes
and relativistic stable processes in $C^{1, 1}$ open
subsets of $\R^d$.
In \cite{BGR},
a Varopoulos type factorization estimate
 in terms of surviving probabilities was obtained for
the transition densities of symmetric stable processes
in $\kappa$-fat open sets.
Very recently, we obtained, in \cite{CKS9},
a Varopoulos type factorization estimate
for the Dirichlet heat kernels
in non-smooth open sets
for a large class of
purely discontinuous subordinate Brownian motions.
We have also obtained in \cite{CKS9} explicit sharp two-sided Dirichlet
heat kernel estimates for a large class of subordinate Brownian motions in
$C^{1,1}$ open sets.

Things become more complicated when
one deals with L\'evy processes
having both Gaussian and jump components.
In \cite{CKS5}, sharped two-sided heat kernel estimates in $C^{1,1}$ open sets
are established for L\'evy processes that can be written as the independent
sum of a Brownian motion and a symmetric $\alpha$-stable process.
A key ingredient is the boundary Harnack principle
for $\Delta +\Delta^{\alpha/2}$
in $C^{1,1}$
open sets with explicit boundary decay rates, obtained in \cite{CKSV}.

The purpose of this paper is to establish sharp two-sided Dirichlet heat kernel
estimates, in $C^{1,1}$ open sets, for a large class of subordinate Brownian motions with Gaussian components.
Throughout this paper,
we will always assume that  $S=(S_t:\ t\ge 0)$
is a complete subordinator with a positive drift and, without loss
of generality, we shall assume that the drift of $S$ is equal to 1.
That is, the Laplace exponent of $S$ is a complete Bernstein function which can be written as
\begin{equation}\label{H:1aa}
\phi(\lambda):=\lambda+\psi(\lambda)
\qquad
\text{ with }
~~
\psi(\lambda):=\int_{(0,\infty)}(1-e^{-\lambda t})\, \mu(dt),
\end{equation}
where $\mu$ is a measure on $(0, \infty)$ satisfying $\int^\infty_0
(1\wedge t)\mu(dt)<\infty$.
$\mu$ is  called the L\'evy measure of the subordinator $S$ (or of $\phi$).
We will exclude the trivial case of $S_t=t$,
 that is,
 the case of $\psi\equiv 0$.
By the definition of complete Bernstein functions,
the L\'evy measure $\mu$ has a complete monotone density.
By a slight abuse of
notation we will denote the density by $\mu(t)$.
For basic facts on complete Bernstein functions, we refer the reader to \cite{SSV}.
In this paper, we will assume the following growth condition on
$\mu(t)$ near zero:
For any $K>0$, there exists $c=c(K)>1$ such that
\begin{equation}\label{H:1a}
\mu(r)\le c \mu(2r), \quad r\in (0, K).
\end{equation}

Suppose that $B=(B_t: t\ge 0)$ is a Brownian motion in $\R^d$ with infinitesimal
generator $\Delta$ and independent of $S$. Then the process $X=(X_t: t\ge 0)$ defined by
$X_t=B_{S_t}$ is called a subordinate Brownian motion. $X$ can be written as the independent
sum of a Brownian motion and a purely discontinuous subordinate Brownian motion.
The infinitesimal generator of $X$ is
$$
\sL^X:= - \phi (-\Delta)=\Delta -\psi (-\Delta),
$$
and $\psi (-\Delta)$ is a non-local operator.
The L\'evy density $J$ of $X$ is given by
\begin{equation}\label{e:levymeasurenew}
J(x)=j(|x|)=\int^\infty_0(4\pi t)^{-d/2}e^{-|x|^2/4t}\mu(t)dt.
\end{equation}
The  function  $J(x)$  determines a L\'evy system for $X$,
which describes the jumps of the process $X$: for any non-negative
measurable function $f$ on $\bR_+ \times \bR^d\times \bR^d$ with
$f(s, y, y)=0$ for all $y\in \bR^d$, any stopping time $T$ (with
respect to the filtration of $X$) and any $x\in \bR^d$,
\begin{equation}\label{e:levy}
\E_x \left[\sum_{s\le T} f(s,X_{s-}, X_s)
\right]= \E_x \left[
\int_0^T \left( \int_{\bR^d} f(s,X_s, y) J(X_s-y) dy \right)
ds \right]
\end{equation}
(see, for example, \cite[Proof of Lemma 4.7]{CK} and \cite[Appendix
A]{CK2}).

The function $j$ is obviously a decreasing function on $(0, \infty)$. Using this and the fact that
$J$ is a L\'evy density (and so $\int_{\R^d}(1\wedge |x|^2)j(|x|)dx<\infty$),
we can easily get that, for any $K>0$, there exists $c=c(K)>0$ such that
\begin{equation}\label{e:jnearorigin}
j(r)\le cr^{-d-2}  \qquad \hbox{for } r\in (0, K].
\end{equation}
In fact, we have for $s\in (0, K]$,
$$
\frac1{d+2}j(s)s^{d+2}=\int^s_0r^{d+1}j(s)dr\le \int^s_0r^{d+1}j(r)dr\le \int^K_0r^{d+1}j(r)dr=:
c(K)<\infty,
$$
from which \eqref{e:jnearorigin} follows immediately.

The subordinate Brownian motion $X$ has a transition density $p(t, x, y)$
with respect to the Lebesgue measure.
Observe that $p(t,x,y)$ is given by $p(t,x,y)=p(t, |x-y|)$ where
\begin{equation}\label{newop}
p(t, r)=\int^\infty_0(4\pi s)^{-d/2} e^{-\frac{r^2}{4s}}\P(S_t\in ds) \quad \forall t>0, r \ge0.
\end{equation}
Clearly $r \to p(t, r)$ is monotonically deceasing.
For any open set $D\subset \R^d$, we will use $X^D$ to denote the part process
of $X$ killed upon leaving $D$.
The process $X^D$ has a transition density $p_D(t, x, y)$  with respect to the Lebesgue measure on $D$.
The density $p_D(t, x,y)$ is the fundamental solution
of $\sL^X$ in $D$ with zero exterior condition, which is also
called the Dirichlet heat kernel of $\sL^X$ in $D$.

The goal of this paper is to establish
explicit sharp two-sided
estimates for $p_D(t, x, y)$ in $C^{1, 1}$ open sets $D$ under the above
assumptions.
Throughout the remainder of this paper, we assume that $d \ge 1 $.
The Euclidean distance between $x$ and $y$
will be denoted as $|x-y|$. We will use $B(x, r)$ to denote the open
ball centered at $x\in \bR^d$ with radius $r>0$.
For $a, b\in \bR$, $a\wedge b:=\min \{a,
b\}$ and $a\vee b:=\max\{a, b\}$.
We also define $J(x, y)=J(y-x)$.
For any Borel set $A\subset \bR^d$, we will use ${\rm diam}(A)$ to denote
its diameter and $|A|$ to denote its
Lebesgue measure.
For any two positive functions $f$ and $g$,
$f\asymp g$ means that there is a positive constant $c\geq 1$
so that $c^{-1}\, g \leq f \leq c\, g$ on their common domain of
definition.

We will show in this paper (see Remark \ref{r:sgloballb} below)
that under the above assumptions,
for every $T>0$, there exists a constant
$c_0=c_0(T, d, \psi) >0$
so that
$$
p(t, r) \geq
  c_0\left(t^{-d/2}e^{-r^2/4t}+ t^{-d/2}\wedge (tj(r))\right)
$$
for all $(t,r)\in (0, T]\times [0, \infty)$.
To get an explicit Dirichlet heat kernel upper bound estimates,
we will need to
assume the following upper bound condition on $p(t, r)$
for $r \le $diam$(D)$:
For any $T>0$, there exist $C_j
\ge 1$, $j=1, 2, 3, $ such that
for all $(t,r)\in (0, T]\times [0, $ diam$(D)]$,
\begin{equation}\label{globalup}
p(t, r) \le C_1 \left(t^{-d/2}e^{-r^2/C_2t}+ t^{-d/2}\wedge (tj(r/C_3))\right).
\end{equation}

In \cite{CK3}, a DeGiorgi-Nash-Moser-Aronson type  theory has been established
for a large class of symmetric
Markov processes on $\R^d$ with infinitesimal generators of the form
\begin{equation}\label{e:1.2}
  \sL  u(x) =  \sum_{i, j=1}^n \frac{\partial }{\partial x_i} \left(a_{ij}(x)
 \frac{\partial u(x)}{\partial x_j}\right) +
 \lim_{\eps \downarrow 0} \int_{\{y\in \R^d: \, |y-x|>\eps\}}
 (u(y)-u(x)) \frac{c(x, y)}{|x-y|^d \,
 \Phi (|x-y|)} dy,
 \end{equation}
where
$(a_{ij}(x))_{1\leq i, j\leq d}$ is a measurable $d\times d$
matrix-valued measurable function on $\R^d$ that is uniformly
elliptic and bounded,  $c(x, y)$ is a measurable
symmetric kernel that is bounded between two positive constants,
and $\Phi (r)$ is a positive increasing function in $r\in (0, \infty)$.
If $\Phi$ satisfies suitable growth conditions near zero and infinity,
sharp two-sided estimates on the transition density of this class of Markov processes
have been obtained in \cite{CK3}.
In this case, the transition density $p(t, x,y)$ of
such a process
admits the following estimates: for any $T>0$, there exist $c_j\ge 1$ such that
for all $(t, x, y)\in (0, T]\times \R^d\times \R^d$,
\begin{eqnarray*}
p(t, x, y)&\ge& c_1^{-1} \left(t^{-d/2}e^{-c_2|x-y|^2/t}+ t^{-d/2}\wedge (tj(c_3|x-y|))\right)\\
p(t, x, y) &\le& c_1 \left(t^{-d/2}e^{-|x-y|^2/c_2t}+ t^{-d/2}\wedge (tj(|x-y|/c_3))\right)
\end{eqnarray*}
with $j(r)=r^{-d} \, (\Phi (r))^{-1}$.
These estimates can be regarded as the counterpart of Aronson's estimates
for non-local operators.
 When $(a_{ij})$ is a constant matrix and
$c(x, y)=c(|x-y|)$ is a function of $|x-y|$,
the Markov process $X$ with generator $\sL$ of \eqref{e:1.2}
is a rotationally symmetric L\'evy process on $\R^d$ with L\'evy measure
$j(|z|)dz
=\frac{c(|z|)}{|z|^d \, \Phi (|z|)}\, dz$.  In this case, the L\'evy exponent of
 $X$ is
\begin{equation}\label{e:1.5a}
 \Psi (\xi )= \sum_{i, j=1}^d a_{ij} \xi_i \xi_j + \int_{\R^d} \left(1-\cos (\xi \cdot z)\right)
   \frac{c(|z|)}{|z|^d \, \Phi (|z|)}\, dz, \qquad \xi \in \R^d .
\end{equation}

For the subordinate Brownian motion $X$ considered in this paper,
its L\'evy exponent is $\phi (|\xi|^2) $, where $\phi$ is defined in \eqref{H:1aa},
which admits an expression of the form
\begin{equation}\label{e:levyexpnew}
\phi (|\xi|^2)=|\xi|^2 + \int_{\R^d} \left(1-\cos (\xi \cdot z)\right)
J(z)dz =|\xi|^2 + \psi(|\xi|^2),
\end{equation}
where $J$ is defined in \eqref{e:levymeasurenew}.
When the complete Bernstein function $\psi$  satisfies the following condition near infinity:
there exist constants $ \delta_k\in (0,1)$, $a_k>0$, $k=1, 2$, and $R_1>0$ such that
\begin{equation}\label{e:1.5}
a_1\lambda^{\delta_1} \psi (r) \leq \psi(\lambda r) \le a_2 \lambda^{\delta_2} \psi (r) \qquad
\hbox{for } \lambda \ge 1 \hbox{ and } r\geq R_1,
\end{equation}
then (see \cite[Lemma 3.2]{KSV5})
$$
J(x)\asymp \frac1{|x|^d\Phi(|x|)}  \qquad \hbox{for } |x|\le 1,
$$
with a positive increasing function $\Phi$
satisfying the growth conditions in \cite{CK3} for $r\leq 1$.
Then one can use the heat kernel estimates in \cite{CK3} and an
argument similar to the proof of \cite[Theorem 2.4]{CKK2} to show that
the  estimate \eqref{globalup} holds for $X$ for any bounded open set $D$
with $C_3=1$.

If,  in addition to \eqref{e:1.5}, $\psi$ also satisfies the following
condition near zero: there exist constants $ \delta_k\in (0,1)$, $a_k>0$, $k=3, 4$, and $R_2>0$ such that
\begin{equation}\label{e:1.4}
a_3\lambda^{\delta_3} \psi (r) \leq \psi(\lambda r) \le a_4 \lambda^{\delta_4} \psi (r) \qquad
\hbox{for } \lambda \le 1 \hbox{ and } r\leq R_2,
\end{equation}
then (see \cite[Theorem 3.4]{KSV5})
$$
J(x)\asymp \frac1{|x|^d\Phi(|x|)}  \qquad \hbox{for } x\neq 0
$$
with a positive increasing function $\Phi$ satisfying
the conditions in \cite{CK3}
for all $r>0$. So
it follows from the heat kernel estimates in \cite{CK3} that
the  estimate \eqref{globalup} holds for $X$ with
$D=\R^d$ and $C_3=1$.

To state the main result of this paper,  we first recall that an open
set $D$ in $\bR^d$ (when $d\ge 2$) is said to be a (uniform)
$C^{1,1}$ open set if there exist a localization radius $R_0>0$ and
a constant $\Lambda_0>0$ such that for every $z\in\partial D$, there
exist a $C^{1,1}$-function $
\varphi=
\varphi_z: \bR^{d-1}\to \bR$
satisfying $
\varphi (0)= 0$, $\nabla
\varphi (0)=(0, \dots, 0)$, $\| \nabla
\varphi  \|_\infty
\leq \Lambda_0$, $| \nabla
\varphi (x)-\nabla
\varphi (z)| \leq \Lambda_0
|x-z|$, and an orthonormal coordinate system $CS_z$ with its origin
at $z$ such that
$$
B(z, R_0)\cap D=\{ y= (\wt y, \, y_d) \mbox{ in } CS_z: |y|< R_0,
y_d
>
\varphi (\wt y) \}.
$$
The pair $(R_0, \Lambda_0)$ is called the characteristics of the
$C^{1,1}$ open set $D$. Note that a $C^{1,1}$ open set $D$ with
characteristics $(R_0, \Lambda_0)$ can be unbounded and
disconnected; the distance between two distinct components of $D$ is
at least $R_0$. Let $\delta_{\partial D}(x)$ be  the Euclidean
distance between $x$ and $\partial D$.
 It is well-known that any $C^{1, 1}$ open set $D$
satisfies both the {\it uniform interior ball condition} and the
{\it uniform exterior ball condition}: there exists
$r_0=r_0(R_0, \Lambda_0)\in (0, R_0]$
 such that for any $x\in D$ with $\delta_{\partial D}(x)< r_0$ and $y\in
\bR^d \setminus \overline D$ with $\delta_{\partial D}(y)<r_0$,
there are $z_x, z_y\in \partial D$ so that $|x-z_x|=\delta_{\partial
D}(x)$, $|y-z_y|=\delta_{\partial D}(y)$ and that $B(x_0,
r_0)\subset D$ and $B(y_0, r_0)\subset \bR^d \setminus \overline D$
for $x_0=z_x+r_0(x-z_x)/|x-z_x|$ and $y_0=z_y+r_0(y-z_y)/|y-z_y|$.
By a $C^{1,1}$ open set in $\bR$ we mean an open set which can be
written as the union of disjoint intervals so that the minimum of
the lengths of all these intervals is positive and the minimum of
the distances between these intervals is positive.

For an open set $D\subset \bR^d$ and $x\in D$, we will use
$\delta_D(x)$ to denote the Euclidean distance between $x$ and
$D^c$.
For an open set $D\subset \bR^d$ and $\lambda_0\in [1, \infty)$, we say {\it the path distance in
each connected component of  $D$ is comparable to the Euclidean
distance with characteristic $\lambda_0$} if for every $x, y$ in the same
component of $D$
there
is a rectifiable curve $l$ in $D$ connecting $x$ to $y$ such that the
length of $l$ is no larger than $\lambda_0|x-y|$. Clearly, such a
property holds for all bounded $C^{1,1}$ open sets, $C^{1,1}$ open
sets with compact complements and
connected open sets above graphs of $C^{1,1}$
functions.

For any open set $D\subset \R^d$ and positive constants $c_1$ and $c_2$, we define
\begin{equation}\label{eq:qd}
h_{D, c_1, c_2}(t, x, y):=
 \left(1\wedge \frac{\delta_D(x)}{\sqrt{t}}\right)
\left(1\wedge \frac{\delta_D(
y)}{\sqrt{t}}\right)
\left(t^{-d/2}e^{-c_1|x-y|^2/t}+ t^{-d/2}\wedge (tJ(c_2x, c_2y))\right).
\end{equation}
The following is the main result of this paper.

\begin{thm}\label{t:main}
Suppose that $X$ is a
subordinate Brownian motion
with L\'evy exponent
$\phi(|\xi|^2)$ with $\phi$ being a complete Bernstein function satisfying
\eqref{H:1aa} and \eqref{H:1a}.
Suppose that $D$ is a $C^{1,1}$
open subset of $\bR^d$  with characteristics $(R_0, \Lambda_0)$.
\begin{description}
\item{\rm (i)}
If the path distance in  each connected component of $D$
is comparable to the Euclidean distance with characteristic $\lambda_0$,
then for every $T>0$, there exist
  $c_1=c_1(R_0, \Lambda_0,   \lambda_0, T,
  \psi, d)>0$ and  $c_2=c_2(R_0, \Lambda_0,   \lambda_0,  d)>0$
such that for all $(t, x, y) \in (0, T]\times D\times D$,
$$
p_D(t, x, y)\,\ge\, c_1 \, h_{D, c_2, 1}(t, x, y).
$$
\item{\rm (ii)}
If $D$ satisfies \eqref{globalup}, then for every $T>0$,
there exists $c_3=
c_3(R_0, \Lambda_0,  T, d,  \psi, C_1, C_3, d )>1$ such
that for all $(t, x, y) \in (0, T]\times D\times D$,
$$
p_D(t, x, y)\,\le\, c_3\,
h_{D, C_4, C_5}(t, x, y),
$$
where  $C_4=(16C_2)^{-1}$ and
$C_5= (8 \vee 4C_3)^{-1}$.
\item{\rm (iii)} If $D$ is bounded,
then for every $T>0$, there exists
$c_4=c_4 ( \text{\rm diam}(D),  R_0, \Lambda_0,
T, \psi, d )>0$
such that for all $(t, x, y)\in [T, \infty)\times
D\times D$,
$$
p_D (t, x, y) \,\ge\,
c_4 \, e^{-\lambda_1 t}\delta_D (x)\delta_D (y),
$$
where $-\lambda_1<0$ is the largest eigenvalue of
the generator of $X^D$.
\item{(iv)} If $D$ is bounded and satisfies \eqref{globalup},
then for every $T>0$, there exists
$c_5=c_5(\text{\rm diam}(D)$, $R_0, \Lambda_0,
 T,  \psi, d, C_1, C_2,$ $C_3)>0$
such that for all $(t, x, y)\in [T, \infty)\times D\times D$,
$$
p_D (t, x, y)
 \,\leq\,
 c_5\, e^{-\lambda_1 t}\, \delta_D (x)\delta_D (y).
$$
\end{description}
\end{thm}

When $D=B(x_0, r)$,
it follows as a special case of \cite[Theorem 4.5(ii)]{CS7}
that $ \phi (\lambda_1^D)/2\leq \lambda_1 \leq  \phi (\lambda_1^D)$,
where $\lambda_1^D$ is the smallest eigenvalue of $-\Delta$ in $D$.
It follows from the scaling property of Brownian motion
 (or Laplacian) that $\lambda_1^{B(x_0, r)}= cr^{-2}$,
where $c=c(d)$ is a positive constant that depends only on the dimension $d$.
When $D$ is a bounded $C^{1,1}$ open set in $\R^d$
 with $C^{1,1}$-characteristics $(R_0, \Lambda_0)$,
 $D$ contains a ball of radius $r_0$ and is contained
in a ball of radius ${\rm diam}(D)$, where $r_0=r_0(R_0, \Lambda_0)$ is
such that $D$ satisfies the uniform interior ball condition with radius $r_0$.
 By the domain monotonicity of
the first eigenvalue $\lambda_1$,  one concludes from above that
$\lambda_1$
is bounded between two positive constants that depend only on
$R_0, \Lambda_0, \psi$, ${\rm diam}(D)$ and $d$.

Note that
the explicit upper bound estimates
in Theorem \ref{t:main} are established under the
assumption that the upper bound \eqref{globalup} for $p(t, x, y)$ holds.
If, instead of \eqref{globalup}, we assume that there exist constants
$\delta\in (0, 1)$ and $C_6>0$ such that
the function $\psi$ in \eqref{H:1aa} has the property
\begin{equation}\label{psiup}
\psi(\lambda r)\le C_6\lambda^\delta\psi(r) \qquad \mbox{for } \lambda\ge 1
\hbox{ and } r\ge 1,
\end{equation}
we can establish
the following upper bound on the Dirichlet heat kernel
in terms of $p(t, x, y)$ and the boundary decay terms.

\begin{thm}\label{t:cks2}
Suppose that $X$ is a subordinate Brownian motion with L\'evy exponent $\phi(|\xi|^2)$ with $\phi$ being a complete Bernstein function satisfying
\eqref{H:1aa}, \eqref{H:1a} and \eqref{psiup}.
Suppose that $D$ is a $C^{1,1}$
open subset of $\bR^d$  with characteristics $(R_0, \Lambda_0)$.
\begin{description}
\item{\rm (i)}
For every $T>0$, there exists
 $c_1=c_1(R_0, \Lambda_0, T,
 \psi, d)>0$
 such that for all $t\in
(0, T]$ and all $x, y\in D$,
$$
p_D(t, x, y) \leq c_1 \,
\left(1\wedge \frac{\delta_D(x)}{\sqrt{t}}\right)\left(1\wedge \frac{\delta_D(y)}{\sqrt{t}}\right)
p(t,x/4,y/4).
$$
\item{\rm (ii)}
If $D$ is bounded,
then for every $T>0$, there exists
$c_2=c_2(\text{\rm diam}(D)$, $R_0, \Lambda_0,
d, T,  \psi) \ge 1$
such that for all $(t, x, y)\in [T, \infty)\times D\times D$,
\begin{equation}\label{e:nlu}
 c_2^{-1}\, e^{-\lambda_1 t}\, \delta_D (x)\delta_D (y) \,\leq\,
p_D (t, x, y)
 \,\leq\,
 c_2\, e^{-\lambda_1 t}\, \delta_D (x)\delta_D (y),
\end{equation}
where $-\lambda_1<0$ is the largest eigenvalue of
the generator of $X^D$.
\end{description}
\end{thm}

By integrating the two-sided heat kernel estimates in
Theorem \ref{t:main} with respect to $t$, we can easily obtain
sharp two-sided estimates on the Green function $ G_D(x, y):=\int_0^\infty
p_D(t, x, y)dt$.
For this, let
\begin{equation}
g_D(x,y) :=\begin{cases} \frac{1} {|x-y|^{d-2}}
\left(1\wedge \frac{  \delta_D(x) \delta_D(y)}{ |x-y|^{2}}
\right)  \qquad &\hbox{when } d\ge3,  \\
\log \left( 1+ \frac{  \delta_D(x) \delta_D
(y)}{ |x-y|^{2}
}\right)  &\hbox{when } d=2, \\
\big( \delta_D(x)  \delta_D (y)\big)^{1/2} \wedge \frac{
\delta_D(x) \delta_D (y)}{ |x-y|} &\hbox{when
} d=1.
\end{cases} \label{e:gD}
\end{equation}

\begin{cor}\label{C:1.5}
Suppose that $X$ is a
subordinate Brownian motion
 with L\'evy exponent $\phi(|\xi|^2)$ with $\phi$ being a complete Bernstein function satisfying
\eqref{H:1aa} and \eqref{H:1a}.
Suppose that $D$ is a bounded $C^{1,1}$
open subset of $\bR^d$
with characteristics $(R_0, \Lambda_0)$.
\begin{description}
\item{\rm (i)}
There
exists $c_1=c_1 ( \text{\rm diam}(D), R_0, \Lambda_0,
   \psi,   d)>0$ such that
$$
G_D(x, y)\,\ge\, c_1 \, g_D(x,y), \qquad x, y\in D.
$$
\item{\rm (ii)}
If $D$ satisfies \eqref{globalup}, then there exists
 $c_2=c_2 ( \text{\rm diam}(D), R_0, \Lambda_0,
\psi, C_1, C_2, C_3, d)>0 $ such that
$$
G_D(x, y)\, \le \, c_2 \, g_D(x,y), \qquad x, y\in D.
$$
\end{description}
\end{cor}
We remark that even though $D$ may be
disconnected, in contrast with the Brownian motion case,
the process $X^{D}$ is always irreducible because $X^{D}$
can jump from one component of $D$ to
another.
Denote by $G^0_D(x,y)$ the Green function
of Brownian motion in $D$.
It is known (see \cite{CZ})
that $G^0_D(x, y) \asymp g_D(x, y)$ when $x$ and $y$ are in the same
component of $D$, and $G^0_D(x, y)=0$ otherwise.
Thus when $D$ is a bounded $C^{1,1}$ connected open subset of $\bR^d$,
the estimates in  Corollary  \ref{C:1.5} are exactly the same as those
for Brownian motion, while our heat kernel estimates (Theorem \ref{t:main}(i)--(ii))
detect  a short-time and short-distance region,  precisely $t\le |x-y|^2\le 1$ and
$\delta_D(x) \wedge \delta_D(y) \ge \sqrt{t} $, where the jump part is the dominant term.
When $\phi (\lambda ) =\lambda + \lambda^{\alpha/2}$,
Theorem \ref{t:main} and Corollary \ref{C:1.5} in particular recover
 the main results of \cite{CKS5}.

Throughout this paper the constants
$r_0$, $R_0$, $\lambda_0$, $\Lambda_0$,
and $C_i$,  $i=1,\dots,
6$, will be fixed.
We use $c_1, c_2, \cdots$ to
denote generic constants, whose exact values are not important and
can  change from one appearance to another. The labeling of the
constants $c_1, c_2, \cdots$ starts anew in the statement of each
result.
We use $c(\alpha, \beta, ...)$ to indicate
a positive constant that depends on the parameters $\alpha$, $\beta$, ...
Dependence on dimension $d$ will not be explicitly mentioned.
We will use $dx$ to denote the Lebesgue measure in
$\bR^d$.

\section{Lower bound estimate}

In this section we derive the lower bound estimate on $p_D(t, x, y)$
when $D$ is a $C^{1,1}$ open set such that
the path distance in  each connected component of $D$
is comparable to the Euclidean distance. As a consequence,
we also get the lower bound estimate on $p(t, x, y)$ in $\R^d$.
We will use some relation between
killed subordinate Brownian motions and subordinate killed Brownian motions.
In this paper we always assume that
$X$ is a
subordinate Brownian motion
with L\'evy exponent
$\phi(|\xi|^2)$ with $\phi$ being a complete Bernstein function satisfying
\eqref{H:1aa} and \eqref{H:1a}.

Let $\widetilde{S}_t$ be a subordinator whose Laplace exponent $\psi$ is given by \eqref{H:1aa}.
Then $t+\widetilde{S}_t$ is a subordinator which has the same law as $S_t$.
Assume that $\widetilde{S}_t$ is independent of the Brownian motion $B$ in $\R^d$.
Suppose that $U$ is an open subset of $\bR^d$.
We denote by $B^U$ the part process of $B$ killed upon leaving $U$.
The process $\{Z^U_t: t\ge 0\}$ defined by $Z^U_t=B^U_{t+\widetilde{S}_t}$
is called a subordinate killed Brownian motion in $U$.
Let $q_U(t, x, y)$ be the transition density of $Z^U$.
Denote by $\zeta^{Z,U}$ the lifetime of $Z^U$. Clearly,  $Z^U_t =B_{t+\widetilde{S}_t}$ for every $t\in [0, \zeta^{Z, U})$.
Therefore we have
\begin{equation}\label{e:pq}
p_{U}(t, z, w)\ge q_{U}(t, z, w) \quad \hbox{for } (t, z, w)\in (0,
\infty)\times U\times U.
\end{equation}
By \cite[Proposition III.8]{Ber},
for every $b>0$,
there exist constants $T_0>0$
and $c>0$ so that
$$
\P(\widetilde{S}_t\le b\, t)>c  \quad \hbox{for } t\le T_0 .
$$
 Using the Markov property of $\wt S_t$, we can easily deduce that
for every $b>0$ and $T>0$, there exists $c=c(b, T, \psi)>0$ such that
\begin{equation}\label{e:2.2}
\P(\widetilde{S}_t\le b \, t)>c \quad \hbox{for } t\le T.
\end{equation}
These facts
(with $b=1$)
 will be used in the proof of the following lemma.

\begin{lemma}\label{lbbyskbm}
Suppose that $D$ is a $C^{1,1}$ open set in $\bR^d$
with characteristics $(R_0, \Lambda_0)$ such that
the path distance in  each connected component of $D$
is comparable to the Euclidean distance with characteristic $\lambda_0$.
For every constant $T>0$, there exist positive constants
$c_1= c_1( R_0, \Lambda_0, \lambda_0, MT, \psi)$ and
$c_2= c_2( R_0, \Lambda_0, \lambda_0)$
such that for all $\lambda \in (0,  M]$, $t\in (0, T]$ and
$x, y$ in the same connected component of $\sqrt{\lambda}D$,
\bee  \label{e:direct}
p_{\sqrt{\lambda}D}(\lambda t, x, y) \ge c_1(\lambda t)^{-d/2} \left(1\wedge \frac{ \delta_{\sqrt{\lambda}D}(x)}
{\sqrt{\lambda t}} \right)\left(1\wedge \frac{ \delta_{\sqrt{\lambda }D}(y)}{\sqrt{\lambda t}} \right)e^{-c_2|x-y|^2/(\lambda t)}.
\eee
\end{lemma}

\pf Suppose that $\lambda^{-1/2} x$ and $\lambda^{-1/2} y$ are in the same component, say $U$, of
$D$. Let $\widetilde{p}_U(t, z, w)$ be the transition density of
$B^U$.
By
\cite[Theorem 3.3]{Ch} (see also \cite[Theorem 1.2]{Zq3}),
there exist positive constants
$c_1=c_1(R_0, \Lambda_0, \lambda_0, T)$ and $c_2=c_2(R_0, \Lambda_0, \lambda_0)$
 such that for any $(s, z, w)\in (0, 2T]
\times U\times
U$,
\bee \label{e:2.41}
\widetilde{p}_U(s, z, w)\ge c_1\left(1\wedge \frac{\delta_U(z)}{\sqrt{s}}
 \right)\left(1\wedge \frac{\delta_U(w)}{\sqrt{s}}
\right)s^{-d/2}e^{-c_2|z-w|^2/s}.
\eee
(Although not
explicitly mentioned in \cite{Ch}, a careful
examination of the proofs in
\cite{Ch} reveals that the constants $c_1$ and $c_2$ in the above lower
bound estimate can be chosen to depend only on
$(R_0, \Lambda_0, \lambda_0, T)$ and $(R_0, \Lambda_0, \lambda_0)$, respectively.)
By using the scaling property of Brownian motion,
we get that, for every $\lambda>0$, $t\in (0, T]$ and $x, y$ in $\sqrt{\lambda}U$,
$$
\wt p_{\sqrt{\lambda}U}(\lambda t, x, y)= \lambda^{-d/2} \wt
p_U (t, \lambda^{-1/2} x, \lambda^{-
1/2}y).
$$
Thus by  \eqref{e:2.41},
\begin{eqnarray}
 \wt p_{\sqrt{\lambda}U}(\lambda t, x, y) &\geq& c_1 (\lambda t)^{-d/2}
\left(1\wedge \frac{\delta_U(\lambda^{-1/2}x)}{\sqrt{t}}
\right)\left(1\wedge \frac{\delta_U(\lambda^{-1/2}y)}{\sqrt{t}} \right)e^{-c_2|x-y|^2/(\lambda t)} \nn\\
&=& c_1 (\lambda t)^{-d/2} \left(1\wedge \frac{ \delta_{\sqrt{\lambda}U}(x)}
{\sqrt{\lambda t}} \right)\left(1\wedge \frac{ \delta_{\sqrt{\lambda }U}(y)}{\sqrt{\lambda t}} \right)e^{-c_2|x-y|^2/(\lambda t)}.\label{e:2.4}
\end{eqnarray}
Now we assume $\lambda \in (0,  M]$.
Recall that $\widetilde{S}_t$ is independent of  $B$ and that $q_{\sqrt{\lambda}U} (t,x,y)$ is the transition  density of
 $Z^{\sqrt{\lambda}U}_t=B^{\sqrt{\lambda}U}_{t+\widetilde{S}_t}$.
Note that for every $0 < t \le T$ and $x, y$ in $\sqrt{\lambda}U$,
$$
q_{\sqrt{\lambda}U}({\lambda t}, x, y)
=\int^{\infty}_{\lambda t}\widetilde{p}_{\sqrt{\lambda}U}(s, x, y)\P({\lambda t}+\widetilde{S}_{\lambda t} \in ds).
$$
So by  \eqref{e:pq}, \eqref{e:2.2} and \eqref{e:2.4},
for every $0 < t \le T$, $\lambda \in (0,  M]$ and $x, y$ in $\sqrt{\lambda}U$,
\begin{eqnarray*}
p_{\sqrt{\lambda}D}({\lambda t}, x, y)
& \ge &
p_{\sqrt{\lambda}U}({\lambda t}, x, y)
 \ge \,  q_{\sqrt{\lambda}U}
({\lambda t}, x, y) \nonumber\\
&\ge& \int^{2{\lambda t}}_{\lambda t}\widetilde{p}_{\sqrt{\lambda}U}(s, x, y)\P({\lambda t}+\widetilde{S}_{\lambda t} \in ds)\nonumber\\
&=&\int^{{\lambda t}}_0\widetilde{p}_{\sqrt{\lambda}U}({\lambda t}+s, x, y)\P(\widetilde{S}_{\lambda t}\in ds)\nonumber\\
& \ge  & c_3
\left(1\wedge \frac{\delta_{\sqrt{\lambda}U}(x)}{\sqrt{t}}  \right)
\left(1\wedge \frac{\delta_{\sqrt{\lambda}U}(y)}{\sqrt{t}}  \right)(\lambda t)^{-d/2} e^{-c_2|x-y|^2/(\lambda t)}
\P(\widetilde{S}_{\lambda t} \le {\lambda t})\nonumber\\
& \ge  & c_4
\left(1\wedge \frac{\delta_{\sqrt{\lambda}U}(x)}{\sqrt{t}}  \right)
\left(1\wedge \frac{\delta_{\sqrt{\lambda}U}(y)}{\sqrt{t}}  \right)(\lambda t)^{-d/2} e^{-c_2|x-y|^2/(\lambda t)}
\nonumber\\
&=&c_4
\left(1\wedge \frac{\delta_{\sqrt{\lambda}D}(x)}{\sqrt{t}}  \right)
\left(1\wedge \frac{\delta_{\sqrt{\lambda}D}(y)}{\sqrt{t}}  \right)(\lambda t)^{-d/2} e^{-c_2|x-y|^2/(\lambda t)}.
\end{eqnarray*}
\qed

\begin{remark}\label{r:globallb}
\rm
Note that the Brownian motion $B$ in $\R^d$ with infinitesimal
generator $\Delta$ has transition density
$$
\wt p(t,x,y)=(4\pi t)^{-d/2} e^{-\frac{|x-y|^2}{4t}}\, ,\quad x,y\in \R^d, \ t>0\, .
$$
Using this instead of \eqref{e:2.41}
and an argument similar to (but easier than) the proof of Lemma \ref{lbbyskbm} with $\lambda=1$
we can get that,
for any  $T>0$, there exists a positive constant
$c= c(T, \psi)$
such that for all $t \in (0, T]$ and
$x, y$ in $\R^d$,
$$
p(t, x, y)\ge c_1 t^{-d/2}
e^{-\frac{|x-y|^2}{4t}}.
$$
\end{remark}

\begin{lemma}\label{l:cks2lemma3.1}
For any positive constants $R$ and $a$, there exists
$c=c(R, a, \psi)>0$ such that for all $z \in \bR^d$ and $r \in (0, R]$,
$$
\inf_{y\in B(z, r)} \P_y \left(
\tau_{B(z, 2r)} > ar^2 \right)\, \ge\,
c.
$$
\end{lemma}

\pf
By Lemma \ref{lbbyskbm},
we have
\begin{eqnarray*}
&& \inf_{y\in B(z, r)} \P_y \left(
\tau_{B(z, 2r)} > ar^2 \right)
 \geq   \P_0 \left( \tau_{B(0, r)} > ar^2 \right)
 = \int_{B(0, r)} p_{B(0, r)} (ar^2, 0, y) dy \\
&\geq & \int_{B(0, r/2)} p_{B(0, r)} (ar^2, 0, y) dy \\
&\geq & \int_{B(0, r/2)} c_1\, (a r^2  )^{-d/2} \left(1\wedge \frac{ \delta_{B(0, r)}(0)}
{\sqrt{a r^2 }} \right)\left(1\wedge \frac{ \delta_{B(0, r)}(y)}{\sqrt{a r^2  }} \right)e^{-c_2|y|^2/(a r^2  )}   dy \\
&=& c_2\, \int_{B(0, 1/2)} a   ^{-d/2} \left(1\wedge \frac{ \delta_{B(0, 1)}(0)}
{\sqrt{a }} \right)\left(1\wedge \frac{ \delta_{B(0, 1)}(z)}{\sqrt{a   }} \right)e^{-c_2|z|^2/a  }   dz \\
&= & c_3(R, a,  \psi )>0 .
\end{eqnarray*}
\qed

Recall that we assume
that \eqref{H:1a} holds.
On the other hand, since $\phi$ is a complete Bernstein function, it follows from
\cite[Lemma 2.1]{KSV3} that there exists $c_1>1$ such that $\mu(t)\le c_1
\mu(t+1)$ for every  $t>1$. Thus by
 \cite[Proposition 13.3.5]{KSV2} and its proof, we have that
for any $K>0$,
there exists $c_2=c_2(K)>1$ such that
\begin{equation}\label{H:1}
j(r)\le c_2\, j(2r), \qquad \forall r\in (0, K],
\end{equation}
and, there exists $c_3>1$ such that
\begin{equation}\label{H:2}
j(r)\le c_3\, j(r+1), \qquad \forall r
\ge 1.
\end{equation}

\begin{lemma}\label{lower bound12}
Suppose that $R>0$ and $b>1$ are constants. Then there  exists
$c = c (R, b, \psi)>0$ such that
for all $r \in (0, R]$, $t\in [r^2/b, br^2]$ and $u, v\in \bR^d$,
\begin{eqnarray*}
p_{B(u, r)\cup B(v,r)}(t , u, v)
\,\ge\, c\,(t^{-d/2}\wedge (tJ(u, v))).
\end{eqnarray*}
\end{lemma}

\pf
Let $r \in (0, R]$, $t\in [r^2/b, br^2]$ and $E= B(u,r)\cup B(v,r)$.
If $|u-v|\le r/2$, by
Lemma \ref{lbbyskbm} (with $T=b$, $\sqrt \lambda =r$ and $D=B(0,1)$) and \eqref{e:jnearorigin},
\begin{align*}
&p_{E}(t, u, v) \ge \inf_{ |z|<r/2}
p_{B(0,r)}(t, 0, z)= \inf_{ |z|<r/2}
p_{B(0,r)}(r^2(t/r^2), 0, z)\\
&\ge c_1 t^{-d/2} \left(1\wedge \frac{ r}
{\sqrt{ t}} \right)\left(1\wedge \frac{ r}{2\sqrt{ t}} \right)e^{-c_2r^2/t}\ge c_3t^{-d/2}\ge c_4(tJ(u, v)\wedge t^{-d/2}).
\end{align*}

If $|u-v|\ge r/2$, we have by the strong Markov property and the
L\'evy system of $X$ in \eqref{e:levy} that
\begin{align*}
p_{E}(t, u, v) &\ge
\E_u
\left[p_{E}(t-\tau_{
B(u,r/8)},
X_{\tau_{
B(u,r/8)}}, v):\tau_{
B(u,r/8)}<t,
X_{\tau_{
B(u,r/8)}}\in
B(v,r/8) \right]\\
&= \int_0^{t} \left(\int_{
B(u,r/8)} p_{
B(u,r/8)}(s, u, w)
\left(\int_{
B(v,r/8)} J(w, z) p_E (t-s, z, v) dz\right) dw \right) ds \\
&\ge\left(\inf_{w\in
B(u,r/8),\, z\in
B(v,r/8)}J(w,z) \right)
\int_0^{t} \P_u \left(\tau_{
B(u,r/8)}>s \right) \left(\int_{
B(v,r/8)}
p_{E}(t-s, z, v) dz \right) ds\\
&\ge  \P_u(\tau_{
B(u,r/8)}>t) \left(\inf_{w\in
B(u,r/8),\, z\in
B(v,r/8)}J(w,z) \right)\int_0^{t}\int_{
B(v,r/8)}p_{
B(v,r/8)}(t-s, z, v)dz ds\\
&= \P_
0(\tau_{
B(0,r/8)}>t)\left( \inf_{w\in
B(u,r/8),\, z\in
B(v,r/8)}
j(|w-z|)\right)
\int_0^{t}\P_
0(\tau_{
B(
0,r/8)} > s)    ds\\
& \ge
t \left(\P_0(\tau_{B(0,r/8)}>br^2)\right)^2
\left(\inf_{w\in B(u,r/8),\,
z\in B(v,r/8)}
j(|w-z|)\right) \\
& \ge  c_5 t  \left(\inf_{
w\in B(u,r/8),\,
z\in B(v,r/8)}
j(|w-z|)\right).\end{align*}
In the last inequality we have used Lemma \ref{l:cks2lemma3.1}.
Note that,
if
$w\in B(u,r/8)$ and
$z\in B(v,r/8)$, then
$$
|w-z| \le |u-w| +|u-v|+|v-z| \le |u-v| + \frac{r}4 \le (2|u-v|) \wedge (|u-v| + \frac{R}4).
$$
Thus using both \eqref{H:1} and \eqref{H:2} we have
\begin{eqnarray*}
p_{E}(t, u, v)
\ge c_6 tj(|u- v|) \ge c_6 (tJ(u, v) \wedge t^{-d/2})\,.
\end{eqnarray*}
The proof is now complete.
\qed

The next lemma in particular implies
that if $x$ and $y$ are in different components of $D$,
the jumping kernel component of the heat kernel
dominates the Gaussian component.

\begin{lemma}\label{L:2.5}
For any given positive constants $c_1, c_2, R$
and $T$, there is a positive constant $c_3=c_3 (R, T, c_1, c_2, \psi)$
so that
 \begin{equation}\label{e:3.22}
t^{-d/2}e^{-r^2/(c_1t)} \leq c_3 (t^{-d/2} \wedge ( t j(c_2 r)) )
\qquad \hbox{for every } r\geq R \hbox{ and } t\in (0, T].
\end{equation}
\end{lemma}

\pf   Observe that \eqref{H:2} implies that there exist $c_4>0$ and $c_5>0$ such that
\begin{equation}\label{e:3.110}
j(c_2 r) \ge c_4 e^{-c_5r}
\quad \hbox{for every } r>1/c_2 .
\end{equation}
For $r> (1/c_2) \vee (2 c_1c_5T)$ and $t\in (0, T]$,
we have $r^2/(2c_1t) > c_5r$ and
\begin{eqnarray*}
t^{-d/2-1}e^{-r^2/(2c_1 t)}
&\le&  t^{-d/2-1}e^{-  ((1/c_2) \vee (2 c_1c_5 T))^2/(2c_1 t)} \\
&\le & \sup_{0<s \le T}s^{-d/2-1}e^{-  ((1/c_2 ) \vee (2 c_1c_5T))^2/(2
c_1 s)}=:c_6 < \infty.
\end{eqnarray*}
 So by \eqref{e:3.110}, when $r> (1/c_2)  \vee (2 c_1c_5T) $ and
 $t \in (0, T]$, we have
\begin{equation}\label{e:3.11}
t^{-d/2}e^{-r^2/(c_1t)} \leq c_6 t e^{-r^2/(2c_1t)} \le  c_6 t e^{-c_5r}
\le (c_6/ c_4)\, t\, j(c_2 r) .
\end{equation}
When $R \le  r \le (1/c_2)  \vee ( 2 c_1c_5T)$ and $t \in (0, T]$,
clearly
 \begin{equation}\label{e:3.111}
t^{-d/2}e^{-r^2/(c_1t)} \leq t
\left(\sup_{s \leq T} s^{-d/2-1}e^{-{R}^2/(c_1s)} \right)
\le c_7 t j(c_2 r).
\end{equation}
The desired inequality \eqref{e:3.22} now follows from \eqref{e:3.11} and \eqref{e:3.111}. \qed

\medskip

Recall that the function $h_{D, c_1, c_2}(t, x, y)$ is
defined in \eqref{eq:qd}.

\begin{thm}\label{t:low}
Suppose that $D$ is a $C^{1,1}$  open set in $\bR^d$
with characteristics $(R_0, \Lambda_0)$ such that the path distance in
each connected component of  $D$
is comparable to the Euclidean distance with characteristic $\lambda_0$.
For every $T>0$, there exist
$c_1=c_1(R_0, \Lambda_0, \lambda_0, T, \psi)>0$ and
$c_2=c_2(R_0, \Lambda_0, \lambda_0)>0$
such that for all $(t, x, y)\in (0, T] \times
D\times D$,
\begin{equation}\label{eq:ppl2}
p_D(t, x, y)\ge c_1h_{D, c_2, 1}(t,x,y).
\end{equation}
\end{thm}

\pf
First note that the distance between two distinct connected components of $D$ is at least
$R_0$.
Hence in view of Lemmas \ref{lbbyskbm} and \ref{L:2.5},
we only need to show that
there exists $c=c(R_0, \Lambda_0, \lambda_0, T, \psi)>0$
such that for all $(t, x, y)\in (0, T] \times
D\times D$,
\begin{equation}\label{eq:ppl22}
p_D(t, x, y)
\ge c\left(\frac{\delta_D(x)}{\sqrt{t}} \wedge 1 \right)
\left(\frac{\delta_D(y)}{\sqrt{t}} \wedge 1 \right) (
(tJ(x, y))
\wedge t^{-d/2}).
\end{equation}

Since $D$ is a $C^{1,1}$ open set,
as mentioned earlier, it satisfies the uniform interior and uniform exterior ball
conditions with radius  $r_0=r_0(R_0, \Lambda_0)\in (0, R_0]$.
Set $T_0=(r_0/4)^2$.
Consequently,  there exists $L=L(r_0)>1$
such that, for all $t \in (0,  T_0]$ and $x,y \in D$, we can choose
$\xi^t_x \in D\cap B(x, L\sqrt t )$ and $\xi^t_y \in
D \cap B(y, L\sqrt t )$ so that
$ B(\xi^t_x, 2\sqrt t)$ and $B(\xi^t_y, 2\sqrt t)$ are subsets of the connected
components  of $D$ that contains $x$ and $y$, respectively.

We first consider the case $t\in (0, T_0]$.
Note that for $u \in B(\xi^t_x,
\sqrt{t})$, we have
$$\delta_D(u) \ge \sqrt{t} \quad \text{and} \quad
|x-u| \le |x-\xi^t_x|+|\xi^t_x-u| \le L \sqrt{t}+\sqrt{t}=(L+1)\sqrt{t}.$$
Thus
by \eqref{e:direct} (with $\lambda=1$), for $t\in (0, T_0]$,
\begin{align}
&\int_{B(
\xi^t_x, \sqrt{t})} p_D({t/3},x,u)du \,\ge\, c_1
\left(\frac{\delta_D(x)}{\sqrt{t}} \wedge 1 \right) \int_{B(
\xi^t_x,
\sqrt{t})} \left(\frac{\delta_D(u)}{\sqrt{t}} \wedge 1 \right)
t^{-d/2}e^{-c_2|x-u|^2/t} du \nonumber\\
&\ge  c_1 \left(\frac{\delta_D(x)}{\sqrt{t}} \wedge 1 \right)
\, t^{-d/2}e^{-c_2 (L+1)^2}\,  |B(
\xi^t_x, \sqrt{t})|\,\ge\, c_3 \left(\frac{\delta_D(x)}{\sqrt{t}} \wedge 1 \right).
\label{e:loww_0}
\end{align}
Similarly, for $t\in (0, T_0]$,
\begin{eqnarray}
\int_{B(
\xi^t_y, \sqrt{t})} p_D({t/3},y,u)du \ge c_3
\left(\frac{\delta_D(y)}{\sqrt{t}} \wedge 1 \right).
\label{e:loww_01}
\end{eqnarray}
By the semigroup property,  for $t\in (0, T_0]$,
\begin{align}
p_D(t, x, y) \geq
\int_{B(
\xi^t_x, {\sqrt{t}})}\int_{B(
\xi^t_y, {\sqrt{t}})}
p_D({t/3}, x, u) p_D({t/3}, u, v)p_D({t/3}, v, y)dudv
\label{e:sg11}.
\end{align}
We consider
the cases $|x-y| \ge \sqrt{t}/8$ and $|x-y| <
\sqrt{t}/8$ separately.

\medskip

\noindent
{\it Case 1}: Suppose $|x-y| \ge \sqrt{t}/8$ and $t\in (0, T_0]$.
Note that by
\eqref{e:sg11}, Lemma~\ref{lower bound12},
symmetry and   \eqref{e:loww_0}--\eqref{e:loww_01},
\begin{align}
&p_D(t, x, y)\nonumber\\
\geq& \int_{B(
\xi^t_y, \sqrt{t})}\int_{B(
\xi^t_x, \sqrt{t})}
p_D(
{t/3}, x, u)p_{B(u, \sqrt{t}/2) \cup
B(v,\sqrt{t}/2)}({t/3}, u, v)p_D
({t/3},v,y)dudv\nonumber\\
\geq& c_4\int_{B(
\xi^t_y, \sqrt{t})}\int_{B(
\xi^t_x, \sqrt{t})}
 p_D({t/3},x ,u)
(t^{-d/2}\wedge (t
J(u,v)))p_D({t/3}, v, y)dudv\nonumber\\
\geq&
c_4 \left(\inf_{(u,v) \in B(
\xi^t_x, \sqrt{t}) \times
B(
\xi^t_y,\sqrt{t})} (t^{-d/2}\wedge (t
J(u,v))) \right)
\int_{B(
\xi^t_y, \sqrt{t})}\int_{B(
\xi^t_x, \sqrt{t})}
p_D(
{t/3},x,u)p_D({t/3},v,y)dudv \nonumber \\
\geq &
 c_4 c_3^2 \left(\inf_{(u,v) \in
B(
\xi^t_x, \sqrt{t}) \times B(
\xi^t_y,\sqrt{t})} (t^{-d/2}\wedge (t
J(u,v))) \right)
\left(\frac{\delta_D(x)}{\sqrt{t}} \wedge 1 \right)
\left(\frac{\delta_D(y)}{\sqrt{t}} \wedge 1 \right). \label{e:loww1}
\end{align}
Since $|x-y| \ge \sqrt{t}/8$,
we have that for $(u,v) \in
B(
\xi^t_x, \sqrt{t}) \times B(
\xi^t_y,\sqrt{t})$,
\begin{eqnarray*}
&&|u-v| \leq |u-
\xi^t_x|+|
\xi^t_x-x|+|x-y|+|y-
\xi^t_y|+|
\xi^t_y-v|\\
&&\leq 2(1+L)\sqrt{t} +|x-y| \le \big(16(1+L)|x-y|\big) \wedge \big(2(1+L)\sqrt{T_0} +|x-y|\big),
\end{eqnarray*}
thus using \eqref{H:1} and \eqref{H:2}
we have
\begin{equation}\label{e:loww2}
\inf_{(u,v) \in B(
\xi^t_x, \sqrt{t}) \times B(
\xi^t_y,\sqrt{t})}
(t^{-d/2}\wedge (t
J(u,v)))\ge
c_5\, (t^{-d/2}\wedge
(tJ(x, y))).
\end{equation}
Thus combining with \eqref{e:loww1} and \eqref{e:loww2},
we conclude that, for $|x-y| \ge \sqrt{t}/8$,
\begin{eqnarray}
p_D(t, x, y)
\ge
c_6\left(\frac{\delta_D(x)}{\sqrt{t}} \wedge 1 \right)
\left(\frac{\delta_D(y)}{\sqrt{t}} \wedge 1 \right) (
(tJ(x, y))
\wedge t^{-d/2}). \label{e:w2}
\end{eqnarray}

\noindent
{\it Case 2}: Suppose $|x-y| <\sqrt{t}/8 $ and $t\in (0, T_0]$.
Then for $(u,v) \in
B(
\xi^t_x, \sqrt{t}) \times B(
\xi^t_y,\sqrt{t})$,
\begin{eqnarray*}
|u-v| \leq  2(1+L)\sqrt{t} +|x-y| \le (2(1+L)+8^{-1})\sqrt{t}.
\end{eqnarray*}
Thus by \eqref{e:direct},
we have for every $(u,v) \in
 B(
\xi^t_x, \sqrt{t}) \times B(
\xi^t_y, \sqrt{t})$,
$$
p_D(t/3,u, v) \ge
c_7\left(\frac{\delta_D(
u)}{\sqrt{t}} \wedge 1 \right)
\left(\frac{\delta_D(
v)}{\sqrt{t}} \wedge 1 \right)t^{-d/2}
e^{-
c_8|u-v|^2/t} \ge
c_9 t^{-d/2}.
$$
Therefore by
\eqref{e:loww_0}--\eqref{e:sg11},
for $t \le T_0$,
\begin{eqnarray}
p_D(t, x, y) &\ge&
c_9c_3^2
\left(\frac{\delta_D(x)}{\sqrt{t}}
\wedge 1 \right)\left(\frac{\delta_D(y)}{\sqrt{t}} \wedge 1 \right)t^{-d/2}\nonumber\\
&\ge&
 c_9c_3^2\left(\frac{\delta_D(x)}{\sqrt{t}} \wedge 1 \right)
\left(\frac{\delta_D(y)}{\sqrt{t}} \wedge 1 \right) \big(
(tJ(x, y))\wedge t^{-d/2}\big).\label{e:w1}
\end{eqnarray}

Combining \eqref{e:w2} and \eqref{e:w1} we get
\eqref{eq:ppl22} for $t \in (0, T_0]$.  When $T>T_0$ and $t\in (T_0, T]$,
observe that
$   T_0/3\leq t-2T_0/3 \leq T-2T_0/3 \leq (T/T_0-2/3)T_0$,
that is, $t-2T_0/3$ is comparable to $T_0/3$ with some universal
constants that depend only on $T$ and $T_0$.
Using  the inequality
\bee \label{e:2.13}
p_D(t, x, y) \geq
\int_{B(
\xi^{T_0}_x, {\sqrt{T_0}})}\int_{B(
\xi^{T_0}_y, {\sqrt{T_0}})}
p_D(T_0/3, x, u) p_D(t-2T_0/3, u, v)p_D(T_0/3, v, y)dudv
\eee
instead of \eqref{e:sg11}
and by considering
the cases $|x-y| \ge \sqrt{T_0}/8$ and $|x-y| <
\sqrt{T_0}/8$ separately, we deduce by the same argument as above
that
\eqref{eq:ppl22}
holds for $t\in [T_0, T]$ and hence for $t\in (0, T]$.
 \qed

\begin{remark}\label{r:sgloballb} \rm
By Lemma \ref{lower bound12}, we have that
for every $T>0$
there is
a positive constant  $c_1=c_1(\psi, T)$
such that for all $(t, x, y)\in (0, T] \times
\R^d\times \R^d$,
\begin{equation}\label{eq:globallbn}
p(t, x, y)\ge p_{B(x, \sqrt t) \cup B(y, \sqrt t)}(t, x, y) \ge c_1\left( t^{-d/2}\wedge (tJ(x, y))\right).
\end{equation}
Together with Remark \ref{r:globallb}, \eqref{eq:globallbn} yields
the following global lower bound on $p(t, x, y)$:  For every $T>0$,
there is
a positive constant  $c_2=c_2(\psi, T)$
such that for all $(t, x, y)\in (0, T] \times
\R^d\times \R^d$,
\begin{equation}\label{eq:globallb}
p(t, x, y)\ge c_2\left(t^{-d/2}e^{-|x-y|^2/(4t)}+ t^{-d/2}\wedge (tJ(x, y))\right).
\end{equation}
\end{remark}

\section{Upper bound estimate}

In this section, we derive the upper bound estimate on $p_D(t, x, y)$
for $C^{1,1}$ open sets satisfying the assumption \eqref{globalup}.
We first record a lemma, Lemma \ref{l:gen1},
which serves as the starting point for the upper bound estimate.
Applying it and using
\eqref{globalup} for $p_D(t, x, y)$
on the right hand side of \eqref{eq:ub},
we can get an intermediate upper bound estimate for
$p_D(t, x, y)$ that has
one boundary decay factor. This is done in Proposition \ref{p:cks1}.
Applying Lemma \ref{l:gen1} again but now using the intermediate upper bound
estimate for $p_D(t, x, y)$ obtained in Proposition \ref{p:cks1}
on
the right hand side of \eqref{eq:ubn}, we can get the desired short time
sharp upper bound estimate for $p_D(t, x, y)$. This is carried out in
the proof of Theorem \ref{t:main}(ii).
Recall that $X$ is a
subordinate Brownian motion
with L\'evy exponent
$\phi(|\xi|^2)$ with $\phi$ being a complete Bernstein function satisfying
\eqref{H:1aa} and \eqref{H:1a}.

\begin{lemma}\label{l:gen1}
Suppose that $U_1,U_3, E$ are open subsets of $\bR^d$
with $U_1, U_3\subset E$ and
${\rm dist}(U_1, U_3)>0$. Let $U_2 :=E\setminus(U_1\cup U_3)$.
If $x\in U_1$ and $y \in U_3$, then for
every $t >0$,
\begin{align}
p_{E}(t, x, y) \,\le\,& \P_x\left(X_{\tau_{U_1}}\in U_2\right)
\left(\sup_{s<t,\, z\in U_2} p_E(s, z, y)\right)
\nn\\&+\int_0^{t} \P_x(\tau_{U_1}>s) \P_y(\tau_E >t-s) ds \left(\sup_{u\in U_1,\, z\in U_3}J(u,z)\right)\label{eq:ubn}\\
\,\le\,& \P_x\left(X_{\tau_{U_1}}\in U_2\right)
\left(\sup_{s<t,\, z\in U_2} p(s, z, y)\right)+ \left(t \wedge \E_x
\left[\tau_{U_1}\right] \right)\left(\sup_{u\in U_1,\, z\in U_3}J(u,z)\right)\label{eq:ub}.
\end{align}
\end{lemma}

\pf
The proof is similar to that of \cite[Lemma 3.4]{CKS5}.
For the reader's convenience,
we spell out the details here.
Using the strong Markov property of $X$, we have
\begin{eqnarray*}
p_{E}(t, x, y) &=&\E_x\left[p_{E}\big(t-\tau_{U_1},
X_{\tau_{U_1}}, y\big)
: \tau_{U_1}<t \right]\\
&=&\E_x\left[p_{E}\big(t-\tau_{U_1},
X_{\tau_{U_1}},
y\big):\tau_{U_1}<t, X_{\tau_{U_1}}\in U_{2} \right] \\
&&~~+ \E_x\left[p_{E}\big(t-\tau_{U_1},
X_{\tau_{U_1}}, y\big):\tau_{U_1}<t,
X_{\tau_{U_1}}\in U_3\right] \,=:\, I\,+\,II\,.
\end{eqnarray*}
Clearly
\begin{eqnarray*}
I \le  \P_x\left( X_{\tau_{U_1}}\in
U_{2}\right) \left( \sup_{s<t,\, z\in U_{2}} p_E(s, z,
y)\right)\le  \P_x\left( X_{\tau_{U_1}}\in
U_{2}\right) \left( \sup_{s<t,\, z\in U_{2}} p(s, z,
y)\right) .
\end{eqnarray*}
On the other hand,
by \eqref{e:levy} and the symmetry,
\begin{eqnarray*}
II&=& \int_0^{t} \left( \int_{U_1} p_{{U_1}}(s, x, u) \left(
\int_{U_3} J(u,z) p_{E}(t-s, z, y) dz\right)
du\right) ds\\
&\le&   \left(\sup_{u\in U_1,\, z\in U_3}J(u,z)\right)
\int_0^{t} \P_x(\tau_{U_1}>s) \left(\int_{U_3}p_{E}(t-s, z,
y)dz
\right) ds\\
&\le&   \left(\sup_{u\in U_1,\, z\in U_3}J(u,z)\right)
\int_0^{t} \P_x(\tau_{U_1}>s) \P_y(\tau_E >t-s) ds.
\end{eqnarray*}
Finally
\begin{eqnarray*}
 \int_0^{t}\P_x(\tau_{U_1}>s) \P_y(\tau_E >t-s)
ds \le\int_0^{t}\P_x(\tau_{U_1}>s) ds \le
t \wedge \E_x [\tau_{U_1}] .
\end{eqnarray*}
This completes the proof of the lemma.
\qed

Recall that $C_1$, $C_2$ and $C_3$ are the constants in \eqref{globalup}.

\begin{prop}\label{p:cks1}
Suppose that
$D$ is a $C^{1,1}$ open set in $\bR^d$ with
characteristics $(R_0, \Lambda_0)$.
Assume that \eqref{globalup} holds.
For every $T>0$, there exists
 $c=c(C_1, C_3, R_0, \Lambda_0, T, \psi)>0$
 such that for all $t\in
(0, T]$ and all $x, y\in D$,
$$
p_D(t, x, y) \leq c \,
\left(1\wedge \frac{\delta_D(x)}{\sqrt{t}}\right)
\left(t^{-d/2}e^{-|x-y|^2/(4C_2t)}+ (t^{-d/2}\wedge tj(|x-y|/(4 \vee 2 C_3)))\right).
$$
\end{prop}

\pf
There exists $r_0=r_0(R_0, \Lambda_0)\in (0, R_0]$ such that $D$
satisfies the uniform interior and uniform exterior ball conditions
with radius $r_0$.
Fix $T>0$ and $t\in (0, T]$. Let $x, y\in D$.
In view of \eqref{globalup}, we only need to show the theorem for
$\delta_{D}(x)<r_0\sqrt{t}/(16\sqrt{T})\le r_0/(16)$,
which we will assume throughout the remainder of this proof.
Choose $x_0\in \partial D$
such that $\delta_{D}(x)=|x-x_0|$.
Let
\begin{align}\label{eq:nn1}
U_1:=B( x_0,
r_0\sqrt{t}/(8\sqrt{T} ) ) \cap D.
\end{align}
Let ${\bf n}(x_0)$ be the unit inward normal of
$D$ at the boundary point $x_0$. Put
\begin{align}\label{eq:nn2}x_1=x_0+\frac{r_0\sqrt{t}}{16\sqrt{T} }{\bf n}(x_0).\end{align}
 Note that
$\delta_{D} (x_1)=\frac{r_0\sqrt{t}}{16\sqrt{T} }$.
Applying the boundary Harnack principle in \cite{KSV4} we get
$$ \P_x(X_{\tau_{U_1}}\in D \setminus U_1) \le c_1
  \P_{x_1}(X_{\tau_{U_1}}\in D \setminus U_1) \frac{\delta_D(x)} {\delta_D(x_1)}
\leq c_1 \frac{16\sqrt{T}\delta_D(x)}{r_0\sqrt{t}}.
$$
Hence
 \begin{equation} \label{e:nnw11}
 \P_x(X_{\tau_{U_1}}\in D \setminus U_1) \leq
   c_2 \left( 1\wedge  \frac{\delta_D(x)}{ \sqrt{t}}\right).
\end{equation}
By \cite[Lemma 4.3]{KSV4},
 \begin{equation}\label{e:nnw2}
\E_x [\tau_{U_1}] \le
 c_3 \sqrt{t}\ \delta_D(x).
 \end{equation}
Thus we have by \eqref{e:nnw11} and \eqref{e:nnw2},
\begin{eqnarray}
\P_{x}\left(\tau_{D}>t/2 \right) &\le& \P_x\left(\tau_{U_1} >t/2 \right) +
\P_x \left( X_{\tau_{U_1}} \in D\setminus U_1\right) \nn\\
&\le&
  \left( \Big(\frac2{t} \E_x\left[\tau_{U_1}\right]\Big) \wedge 1\right)
 +\P_x \left( X_{\tau_{U_1}} \in D\setminus U_1\right)
\le c_4\left(1 \wedge \frac{\delta_D(x)} { \sqrt{t}}\right)\label{e:nnw22} .
\end{eqnarray}

Now we deal with two cases separately.

\noindent
{\it Case 1}:
$|x-y|\le  \left( \sqrt{2dC_2}\vee (r_0/\sqrt{T})\right)
\sqrt{t}$.
By the semigroup property, symmetry and \eqref{globalup},
\begin{eqnarray*}
 p_{D}(t,x,y) &=&\int_{D} p_{D}({t/2},x,z)
p_{D}({t/2},z,y)  dz\\
&\le& \left(\sup_{z,w \in D}
p({t/2},z,w) \right)\int_{D}
p_{D}({t/2},x,z) dz\\
&\le&
2^{1+d/2}C_1 t^{-d/2}\P_{x}(\tau_{D}>{t/2})\nn\\
&\le& c_42^{1+d/2}C_1t^{-d/2} \left(1 \wedge \frac{\delta_D(x)} { \sqrt{t}}\right),
\end{eqnarray*}
where in the last line   \eqref{e:nnw22} is used.
Since
\begin{align}\label{eq:nn4}
|x-y|^2/(4C_2t) \le ( d/2)\vee (r_0^2/(4C_2T)) \le ( d/2)\vee (r_0^2/(4T)),
\end{align}
we have
\begin{eqnarray*}
 p_{D}(t,x,y) \le  c_42^{1+d/2}C_1e^{( d/2)\vee (r_0^2/(4T))}t^{-d/2}e^{-|x-y|^2/(4C_2t)} \left(1 \wedge \frac{\delta_D(x)} { \sqrt{t}}\right).
\end{eqnarray*}

\noindent
{\it Case 2}:
$|x-y|\ge \left( \sqrt{2dC_2}\vee (r_0/\sqrt{T}) \right) \sqrt{t} $.
 Let
\begin{align}\label{eq:nn5}
U_3:= \{z\in D: |z-x|> |x-y|/2\} \quad \text{ and }\quad
 U_2:=D\setminus (U_1\cup U_3).
\end{align}
Since
$$
|z-x| > \frac{|x-y|}2  \ge \frac{r_0\sqrt{t}}{2\sqrt{T}} \quad \text{for
}z\in U_3,
$$
${\rm dist}(U_1, U_3)>0$ and, for $u\in U_1$ and $z\in U_3$,
\begin{eqnarray}\label{e:n001}
|u-z| \ge |z-x|-|x-x_0|-|x_0-u| \ge |z-x|- \frac{r_0\sqrt{t}}{4\sqrt{T}} \ge
 |z-x| /2 \ge  |x-y| /4 .
\end{eqnarray}
Thus,
\begin{eqnarray}
\sup_{u\in U_1,\, z\in U_3}J(u,z) \le c \sup_{(u,z):|u-z| \ge
\frac{1}{4}|x-y|}j(|u-z|)  \, \le \, c_3 j( |x-y|/4)  .
\label{e:n01}
\end{eqnarray}
If $z \in U_2$,
\begin{equation}\label{e:one2}
\frac32 |x-y| \ge |x-y| +|x-z| \ge  |z-y| \ge |x-y| -|x-z| \ge
\frac{|x-y|}2 .
\end{equation}
We remark here that up to this point, we have not
used assumption
\eqref{globalup} yet in this proof.

Observe that, for any
$\beta>0$, the function $f(s):=s^{-d/2}e^{-\beta/s}$ is increasing
 on the interval $(0, 2\beta /d]$.
By \eqref{globalup},
\eqref{e:one2} and the observation that $t  \le |x-y|^2/(2dC_2)$,
\begin{eqnarray}
\sup_{s\le t,\, z\in U_2} p(s, z, y)&\leq & C_1
\sup_{s\le t,\, z\in U_2} \big(s^{-d/2}e^{-|z-y|^2/C_2s}+ s^{-d/2}\wedge s
J(z/C_3, y/C_3)\big)\nn\\
&\le &C_1
\sup_{s\le t,\, |z-y| \ge |x-y|/2} \big(s^{-d/2}e^{-|z-y|^2/C_2s}+
s
J(z/C_3,y/C_3)\big)\nn\\
& \le &
C_1 \sup_{s\le t} s^{-d/2}e^{-|x-y|^2/(4C_2s)} +
C_1 t
j(|x-y|/(2C_3))\nn\\
&\le&
C_1t^{-d/2}e^{-|x-y|^2/(4C_2t)} + c_5 (t^{-d/2}\wedge t
j(|x-y|/(2C_3))), \label{e:n02}
\end{eqnarray}
where in the last line   \eqref{e:jnearorigin} is used. In fact,
since $|x-y|\ge  (r_0/\sqrt{T}) \sqrt{t}$, by \eqref{e:jnearorigin},
\begin{eqnarray}
tj( |x-y|/ (2C_3)) &\le& t j( (|x-y|\wedge r_0)/(2C_3) )
\le c_6 \left( \frac{t}{|x-y|^2 \wedge r_0^2}\right)^{1+d/2} t^{-d/2}
\nn \\
& \le& c_6(T/r_0^2)^{1+d/2}\, t^{-d/2}, \label{e:nnw112}
\end{eqnarray}
where $c_6>0$ depends only on $C_3$.

By the same argument as that used to get \eqref{e:nnw11},
we can apply the boundary Harnack principle in \cite{KSV4} to get
\bee \label{e:nnw1}
 \P_x(X_{\tau_{U_1}}\in U_2)  \le  c_7
  \P_{x_1}(X_{\tau_{U_1}}\in U_2) \frac{\delta_D(x)} {\delta_D(x_1)}
 \le  c_8\frac{\delta_D(x)} {\sqrt{t}}.
\eee
Applying
\eqref{eq:ub},  \eqref{e:nnw2},
\eqref{e:n01},
\eqref{e:n02} and  \eqref{e:nnw1}, we obtain
\begin{align*}
&p_{D}(t, x, y) \nn\\
\le&c_{9} \left( t^{-d/2}e^{-|x-y|^2/(4C_2t)} +  t^{-d/2}\wedge t
j \left( |x-y|/(2C_3) \right)  \right)
\frac{\delta_D(x)}{\sqrt{t}}+ c_{10}\ tj( |x-y|/4)\frac{\delta_D(x)}{\sqrt{t}}\nn\\
\le& \, c_{11} \left(t^{-d/2}e^{-|x-y|^2/(4C_2t)} +  t^{-d/2}\wedge t
j( |x-y|/(2C_3\vee 4) ) \right)
\frac{\delta_D(x)}{\sqrt{t}},
\end{align*}
where in the last line \eqref{e:jnearorigin} is used (see \eqref{e:nnw112}).
This combined with \eqref{globalup} completes the
proof of this proposition.
\qed

\begin{prop}\label{p:cks1n}
Suppose that
$D$ is a $C^{1,1}$ open set in $\bR^d$ with
characteristics $(R_0, \Lambda_0)$.
Assume that \eqref{globalup} holds.
For every $T>0$, there exists
 $c=c(C_1, C_3, R_0, \Lambda_0, T, \psi)>0$
 such that for all $t\in
(0, T]$ and all $x\in D$
$$ \P_x(\tau_D >t)  \leq c \,
\left(1\wedge \frac{\delta_D(x)}{\sqrt{t}}\right).
$$
\end{prop}
\pf
Fix $T>0$.
By Proposition \ref{p:cks1} and
\eqref{eq:globallb} we have that for every  $0 < t \le T$ and $x, z$ in $D$,
\begin{align*}
p_{D}({t},x,z)
\le& c_1
\left(1\wedge \frac{\delta_D(x)}{\sqrt{t}}\right)
\left(t^{-d/2}e^{-|x-z|^2/C_2   (4t)}+ t^{-d/2}\wedge (tj
\left( |x-z|/(4 \vee 2 C_3) \right)\right)\\
\le& c_3 \left(1\wedge \frac{\delta_D(x)}{\sqrt{t}}\right) p(t, c_2x,c_2z),
\end{align*}
where $c_2:=(\sqrt{C_2} \vee 4 \vee (2 C_3))^{-1}$.
Thus
$$ \P_x(\tau_D >t)  =\int_D p_{D}({t},x,z)dz
\le c_3\left(1\wedge \frac{\delta_D(x)}{\sqrt{t}}\right) \int_Dp(t, c_2x,c_2z) dz
 \le  c_3\left(1\wedge \frac{\delta_D(x)}{\sqrt{t}}\right).$$
\qed

\noindent{\bf Proof of Theorem \ref{t:main}:}
\noindent{(i)}.
 This has already been established in Theorem \ref{t:low}.

\medskip

\noindent{(ii)}.
Fix $T>0$.
There exists $r_0=r_0(R_0, \Lambda_0)\in (0, R_0]$ such that $D$
satisfies the uniform interior and uniform exterior ball conditions
with radius $r_0$.  Let $t\in (0, T]$ and $x, y\in D$.
By Proposition \ref{p:cks1}, \eqref{globalup} and symmetry, we only need to prove (ii) for
$\delta_{D}(x) \vee \delta_{D}(y)<r_0\sqrt{t}/(16\sqrt{T})\le r_0/(16)$,
which we will assume throughout the remainder of
the proof of (ii).

The proof of (ii)
 is along the line of the proof of Proposition \ref{p:cks1}
but using the estimate from Proposition \ref{p:cks1}
for the upper bound estimate of $p_D(t, x, y)$ on the right hand side
of \eqref{eq:ubn} rather than using \eqref{globalup}. Define
$U_1$, $x_0$ and $x_1$
in the same way as in the proof of Proposition \ref{p:cks1}
(see \eqref{eq:nn1}--\eqref{eq:nn2}), and consider the following
two cases separately.

\noindent
{\it Case 1}:
$|x-y|\le  \left( \sqrt{8(d+1)C_2}\vee (r_0/\sqrt{T})\right)
\sqrt{t}$.
By the semigroup property, symmetry and Proposition \ref{p:cks1},
\begin{eqnarray}
 p_{D}(t,x,y) &=&\int_{D} p_{D}({t/2},x,z)
p_{D}({t/2},z,y)  dz \nn \\
&\le& \left(\sup_{z \in D} p_D({t/2},y,z) \right)\int_{D}
p_{D}({t/2},x,z) dz \nn \\
&\le&
c_1 t^{-d/2}\left(1 \wedge \frac{\delta_D(y)} { \sqrt{t}}\right)\P_{x}(\tau_{D}>{t/2})\nn\\
&\le& c_1t^{-d/2}\left(1 \wedge \frac{\delta_D(y)} { \sqrt{t}}\right) \left(1 \wedge \frac{\delta_D(x)} { \sqrt{t}}\right)
\nn \\
&\leq &  c_2t^{-d/2}e^{-|x-y|^2/(16C_2t)} \left(1 \wedge \frac{\delta_D(x)} { \sqrt{t}}\right)\left(1 \wedge
\frac{\delta_D(y)} { \sqrt{t}}\right) , \label{eq:nn12}
\end{eqnarray}
where in the third inequality,  Proposition \ref{p:cks1n} is used.

\medskip

\noindent
{\it Case 2}:
$|x-y|\ge \left( \sqrt{8(d+1)C_2}\vee (r_0/\sqrt{T}) \right) \sqrt{t} $.
Define $U_2$ and $U_3$ as in \eqref{eq:nn5}. Note that \eqref{e:n01} holds.
Moreover, for $z \in U_2$, as in \eqref{e:one2},
\begin{equation}\label{e:one2n}
3 |x-y| /2  \ge  |z-y| \ge |x-y|/2 .
\end{equation}
Observe that, for any
$\beta>0$, the function $f(s):=s^{-(d+1)/2}e^{-\beta/s}$ is increasing
 on the interval $(0, 2\beta /(d+1)]$.
By Proposition \ref{p:cks1} (instead of \eqref{globalup}),
\eqref{e:one2n}
 and the observation that $t  \le |x-y|^2/(8(d+1)C_2)$,
\begin{eqnarray}
\sup_{s\le t,\, z\in U_2} p_D(s, z, y)&\leq &
c_3
\sup_{s\le t,\, z\in U_2} \big(s^{-d/2}e^{-|z-y|^2/(4C_2s)}+ s^{-d/2}\wedge s
j(|z-y|/(4\vee 2C_3))\big)\frac{\delta_D(y)} { \sqrt{s}}\nn\\
&\le &c_3\delta_D(y)
\sup_{s\le t,\, |z-y| \ge |x-y|/2} \big(s^{-(d+1)/2}e^{-|z-y|^2/(4C_2s)}+
\sqrt{s}
j(|z-y|/(4\vee 2C_3))\big)\nn\\
& \le &c_3\delta_D(y)
\big( \sup_{s\le t} s^{-(d+1)/2}e^{-|x-y|^2/(16C_2s)} +
\sqrt{t}
j(|x-y|/(8\vee 4C_3))\big)\nn\\
&\le&c_4\frac{\delta_D(y)}{\sqrt{t}}\big(
t^{-d/2}e^{-|x-y|^2/(16C_2t)} + (t^{-d/2}\wedge t
j(|x-y|/(8\vee 4C_3))\big), \label{e:n02n}
\end{eqnarray}
where in the last line we used an argument similar to that in \eqref{e:nnw112}.
On the other hand,
by Proposition \ref{p:cks1n}
we have
\begin{align}
&\int_0^{t} \P_x(\tau_{U_1}>s) \P_y(\tau_D >t-s) ds
\le \int_0^{t} \P_x(\tau_{D}>s) \P_y(\tau_D >t-s) ds \le c_5 \int_0^{t}\frac{\delta_D(x)}{\sqrt{s}} \frac{\delta_D(y)}{\sqrt{t-s}}ds \nn\\
&= c_5 \delta_D(x) \delta_D(y) \int_0^1 \frac1{\sqrt{r(1-r)}} dr
= c_6 \delta_D(x) \delta_D(y)  . \label{eq:nn8}
\end{align}
Combining
\eqref{eq:ubn},
\eqref{e:n01},
\eqref{e:nnw1},  \eqref{e:n02n} and \eqref{eq:nn8} all together, we conclude that
\begin{eqnarray*}
&& p_{D}(t, x, y) \\
&\le& c_{7} \left( t^{-d/2}e^{-|x-y|^2/(16C_2t)} + (t^{-d/2}\wedge t
j( |x-y|/(8\vee 4C_3))) \right)
\frac{\delta_D(x)\delta_D(y)}{t} \\
&& + c_{8}\ tj(|x-y|/4)\frac{\delta_D(x)\delta_D(y)}{t}\\
&\le& \, c_{9} \left(t^{-d/2}e^{-|x-y|^2/(16C_2t)} + (t^{-d/2}\wedge t
j(|x-y|/(8\vee 4C_3)))\right) \frac{\delta_D(x)\delta_D(y)}{t}
\nn\\
&=& \, c_{9} \left(t^{-d/2}e^{-|x-y|^2/(16C_2t)} + (t^{-d/2}\wedge t
j(|x-y|/ (8\vee 4C_3)))\right)
\left(1 \wedge \frac{\delta_D(x)} { \sqrt{t}}\right) \left(1 \wedge \frac{\delta_D(y)} { \sqrt{t}}\right),
\end{eqnarray*}
where in the second inequality \eqref{e:jnearorigin} is used (see \eqref{e:nnw112}).
This combined with \eqref{eq:nn12} and Proposition \ref{p:cks1}  completes  the proof of (ii).

\medskip

\noindent{(iii)} and (iv).
We first note that the path distance condition is satisfied
in any bounded $C^{1,1}$ open set $D$
with $\lambda_0$ depending only on
$R_0$, $\Lambda_0$ and ${\rm diam} (D)$.
Thus by (i) and (ii), it suffices to prove (iii) and (iv) for $T=3$.

In view of \eqref{globalup},
the transition semigroup $\{P^D_t, t>0\}$
of $X^D$
consists of Hilbert-Schmidt operators, and hence compact operators, in $L^2(D; dx)$.
So $P^D_t$ has discrete spectrum
$\{e^{-\lambda_k t}; k\geq 1\}$ arranged in decreasing order and repeated
according to their multiplicity.
Let $\{\phi_k, k\geq 1\}$ be the corresponding eigenfunctions with
unit $L^2$-norm, which forms an orthonormal basis for $L^2(D; dx)$.

 Clearly,
\begin{equation}\label{ne:4.2}
\int_D \left( 1\wedge  \delta_D(x) \right)
\phi_1 (x) dx   \le
|D|^{1/2} \|\phi_1\|_{L^2 (D)} \leq
|D|^{1/2}.
\end{equation}
By using the eigenfunction expansion of $p_D(t, x, y)=\sum_{k=1}^\infty
e^{-\lambda_k t} \phi_k (x) \phi_k (y)$, we get
\begin{equation}\label{ne:4.3}
\int_{D\times D} \left( 1\wedge  \delta_D(x) \right)
 p_D(t, x, y)  \left( 1\wedge  \delta_D(y) \right)
  \, dx dy =
\sum_{k=1}^\infty e^{-t \lambda_k} \left(
\int_D \left( 1\wedge  \delta_D(x) \right)
\phi_k (x) dx\right)^2.
\end{equation}
 Noting that
$\lambda_k$ is increasing and $\| f\|_2^2 = \sum_{k=1}^\infty (\int_D f(z)\phi_k(z)dz)^2$,
 we have
for all $t>0$,
\begin{align}
\int_{D\times D} \left( 1\wedge  \delta_D(x) \right)
 p_D(t, x, y)  \left( 1\wedge  \delta_D(y) \right)
  \, dx dy
&\leq
e^{-t \lambda_1} \, \int_D \left( 1\wedge  \delta_D(x) \right)^2 dx
 \leq e^{-t \lambda_1} \, |D|.\label{ne:4.4}
\end{align}

On the other hand,
by Theorem \ref{t:main}(ii)
and \eqref{ne:4.2},
there exists $
c_1>0 $ so that for every  $x\in D$,
\begin{eqnarray}
\phi_1 (x)&=& e^{\lambda_1  }
\int_D p_D (1, x, y)
\phi_1 (y) dy\nn\\
&\leq&
c_1 e^{\lambda_1  } \left( 1\wedge  \delta_D(x) \right)
\int_D \left( 1\wedge  \delta_D(y) \right)
\phi_1 (y) dy
\,\leq \,
c_1 e^{\lambda_1  } |D|^{1/2}  \, \left( 1\wedge  \delta_D(x) \right) .
\end{eqnarray}
It now follows from \eqref{ne:4.3}
that for every $t>0$,
\begin{eqnarray}\label{ne:4.55}
&&  \int_{D\times D} \left( 1\wedge  \delta_D(x) \right)
 p_D(t, x, y)  \left( 1\wedge  \delta_D(y) \right)
  \, dx dy \nn \\
 &\geq&
e^{-t \lambda_1} \, \left(\int_D  \left( 1\wedge  \delta_D(x) \right) \phi_1(x) dx
\right)^2 \nonumber  \\
&\geq &  e^{-t \lambda_1} \, \left(\int_D
(c_1 e^{\lambda_1  } |D|^{1/2})^{-1}
\phi_1(x)^2 dx\right)^2 =
c_1^{-2} |D|^{-1}
 \, e^{-(t+2) \lambda_1} .
\end{eqnarray}

For $t\geq
3$ and $x, y\in D$, we have  that
\begin{equation}\label{ne:4.6}
p_D (t, x, y)=\int_{D\times D}
p_D (1, x, z) p_D(t-2 , z, w)
 p_D(1, w, y) dz dw .
\end{equation}
By Theorem \ref{t:main}(ii)
and \eqref{ne:4.4}, there exist
$c_i>0$, $i=2, 3$,
so that for every
$t\geq 3$ and $x, y\in D$,
\begin{eqnarray}
&& p_D(t, x, y)\nn\\
 &\leq&
 c_2 \left( 1\wedge  \delta_D(x) \right) \left( 1\wedge  \delta_D(y) \right) \int_{D\times D}
 \left( 1\wedge  \delta_D(z) \right) p_D (t-2, z, w) \left( 1\wedge  \delta_D(w) \right) dz dw \nn\\
  &\leq&  c_2 |D| e^{-\lambda_1 (t-2)} \, \left( 1\wedge  \delta_D(x) \right) \left( 1\wedge  \delta_D(y) \right)
 \,\le   c_3 e^{-\lambda_1 t}
 \left( 1\wedge  \delta_D(x) \right) \left( 1\wedge  \delta_D(y) \right).
  \label{ne:4.7}
\end{eqnarray}
By \eqref{ne:4.6},
Theorem \ref{t:low}, the boundedness of $D$ and \eqref{ne:4.55}  we have
that there exist
$c_i>0$, $i=4, 5$,
so that for every  $t\geq
3$ and $x, y\in D$,
\begin{eqnarray*}
&& p_D(t, x, y) \nn \\
&\geq&
c_4\, \left(1\wedge j(\text{\rm diam}(D) )\right)^2\left( 1\wedge  \delta_D(x) \right) \left( 1\wedge  \delta_D(y) \right) \int_{D\times D}
 \left( 1\wedge  \delta_D(z) \right) p_D (t-2 , z, w)  \left( 1\wedge  \delta_D(w) \right) dz dw \nn\\
& \geq &
c_5 \left(1\wedge j(\text{\rm diam}(D) )\right)^2  \, |D|^{-1} \left( 1\wedge  \delta_D(x) \right) \left( 1\wedge  \delta_D(y) \right) e^{-t\lambda_1}=
c_6 \left( 1\wedge  \delta_D(x) \right) \left( 1\wedge  \delta_D(y) \right) e^{-t\lambda_1}.
\end{eqnarray*}
 The theorem is now proved.    \qed

Let $(\sE, \sF)$ be the Dirichlet form  of $X$
on $L^2(\R^d; dx)$. It is known that $(\sE, \sF)$ is a regular Dirichlet
form on $L^2(\R^d; dx)$ with core $C^1_c(\R^d)$; see \cite{CF}.
Moreover, for $u\in C^1_c(\R^d)$,
 \begin{equation}\label{eqn:DE}
 \sE (u, u):=\int_{\R^d} \nabla u(x)\cdot  \nabla v(x)dx+
\int_{\R^d\times \R^d} (u(x)-u(y))^2 J(x, y) dx dy
\end{equation}
and  $\sF \,:= \, \overline{C^1_c (\R^d)}^{\sE_1}\subset
 W^{1,2}(\R^d) =\{f\in L^2(\R^d; dx): \sE(f,f)<\infty\}$.
So we have the following Nash's inequality:
 \begin{equation}\label{e:Nash1}
\|f\|_2^{2+4/d}\,\le\, c_1\,\int_{\mathbb R^d} |\nabla
u(x)|^2dx\cdot\|f\|_1^{4/d}\,\le\, c_2\,\sE (f,f)\|f\|_1^{4/d}\qquad
\hbox{for }
f\in \sF.
 \end{equation}
It follows then
\begin{equation}\label{globalupn}
p(t, x,y)\leq c_3
  t^{-d /2}   \quad \hbox{for   }  t>0
 \hbox{ and } x, y \in \R^d.
\end{equation}

\bigskip \noindent
{\bf Proof of Theorem \ref{t:cks2}}. (i)
There exists $r_0=r_0(R_0, \Lambda_0)\in
(0, R_0]$
so that $D$
satisfies the uniform interior and uniform exterior ball conditions
with radius $r_0$.
Fix $T>0$. We claim that there is a constant $c_0>0$ so that
\begin{equation}\label{e:3.30}
 p_{D}(t, x, y)
\leq c_0 p(t, |x-y|/4)
\left(\frac{\delta_D(x)}{\sqrt{t}} \wedge 1\right)
\quad \hbox{for every } (t, x, y)\in (0, T]\times D\times D.
\end{equation}
In view of \eqref{newop},
$p_D (t, x, y) \leq p(t, |x-y|)\leq
p (t, |x-y|/4)$.
So it suffices to prove \eqref{e:3.30} when
  $\delta_{D}(x)<r_0\sqrt{t}/(16\sqrt{T})\le r_0/16 $.
Define
$U_1$, $x_0$ and $x_1$
in the same way as in the proof of Proposition \ref{p:cks1}
(see \eqref{eq:nn1}--\eqref{eq:nn2}).
Hence \eqref{e:nnw11}--\eqref{e:nnw22} and \eqref{e:nnw1} hold.
We now prove \eqref{e:3.30} by considering the following two cases.

\noindent
{\it Case 1}:
$|x-y|\le  (r_0/\sqrt{T}) \sqrt{t}$.
By the semigroup property, symmetry and \eqref{globalupn},
\begin{eqnarray*}
 p_{D}(t,x,y) &=&\int_{D} p_{D}({t/2},x,z)
p_{D}({t/2},z,y)  dz\\
&\le& \left(\sup_{z,w \in D}
p({t/2},z,w) \right)\int_{D}
p_{D}({t/2},x,z) dz\\
&\le&
c_1 t^{-d/2}\P_{x}(\tau_{D}>{t/2})\,\le\, c_2t^{-d/2} \left(1 \wedge \frac{\delta_D(x)} { \sqrt{t}}\right),
\end{eqnarray*}
where in the last line   \eqref{e:nnw22} is used.
Since
$
|x-y|^2/(64t) \le r_0^2/(64T),
$
we have by Remark \ref{r:sgloballb},
\begin{eqnarray*}
p_{D}(t,x,y) \le  c_2e^{r_0^2/(64T)}t^{-d/2}e^{-|x-y|^2/(64t)} \left(1 \wedge \frac{\delta_D(x)} { \sqrt{t}}\right) \le c_3 p(t, |x-y|/4)  \left(1 \wedge \frac{\delta_D(x)} { \sqrt{t}}\right).
\end{eqnarray*}

\noindent
{\it Case 2}:
$|x-y|\ge (r_0/\sqrt{T})  \sqrt{t} $.
 Define $U_2$ and $U_3$ as in \eqref{eq:nn5}.
  Note that \eqref{e:n01} and  \eqref{e:one2} hold.
Observe that by \eqref{newop}
\begin{equation}\label{newoa}
\sup_{s\le t,\, z\in U_2} p(s, z,y) \le \sup_{s\le t,\, |z-y| \ge |x-y|/2} p(s, z,y )
\le \sup_{s\le t } p(s, |x-y|/2).
\end{equation}
 By \cite[Lemma 3.1]{KSV5},  \eqref{psiup} implies
$$
    j(r)\le c_2 r^{-d}\psi(r^{-2})  \le c_1c_2 \psi(1) r^{-d-2\delta}  \quad \text{ for all }r \in (0,1].
$$
Thus under assumption \eqref{psiup}, according to \cite[Theorem 1.3]{CK3},
the parabolic Harnack inequality
holds for the subordinate Brownian motion $X$.
 Extend the definition of $p(t, r)$ by setting $p(t, r)=0$ for $t<0$ and $r\geq 0$. For each fixed $x, y\in \R^d$ and $t>0$ with $|x-y|\ge (r_0/\sqrt{T})  \sqrt{t}$, one can easily check that $(s, w)\mapsto p(s,  |w-y|/2)$
is a parabolic function in $(-\infty, T]\times B(x,  (r_0/\sqrt{T})  \sqrt{t}/4  )$.
So by the parabolic Harnack inequality from \cite[Theorem 1.3]{CK3},
 there is a constant
$c_3=c_3(
  \psi)\geq 1$ so that for every $t\in (0, T]$,
$$  \sup_{s\le t } p(s, |x-y|/2) \le c_3
 p(t, |x-y|/2) .
$$
Hence we have
\begin{eqnarray}
\sup_{s\le t,\, z\in U_2} p(s, z, y)
\leq   c_3  p(t, |x-y|/2). \label{e:n02nn}
\end{eqnarray}
Applying
\eqref{eq:ub},  \eqref{e:nnw2},
\eqref{e:n01}, \eqref{e:nnw1},
\eqref{e:n02nn}, we obtain
\begin{eqnarray}
 p_{D}(t, x, y)
&\le&
\P_x\left(X_{\tau_{U_1}}\in U_2\right)
\left(\sup_{s<t,\, z\in U_2} p(s, z, y)\right)+ \E_x
\left[\tau_{U_1}\right]\left(\sup_{u\in U_1,\, z\in U_3}J(u,z)\right)\nn\\
&\le&
c_{5} p(t, |x-y|/2)
\frac{\delta_D(x)}{\sqrt{t}}+ c_{5}\ tj( |x-y|/4)\frac{\delta_D(x)}{\sqrt{t}}.\label{e:nnw1122n}
\end{eqnarray}
Since $|x-y|\ge   r_0 \sqrt{t} /\sqrt{T}  $,
 by \eqref{e:jnearorigin},
\begin{eqnarray}
tj( |x-y|/ 4) &\le& t j( (|x-y|\wedge r_0)/4 )
\le c_4 \left( \frac{t}{|x-y|^2 \wedge r_0^2}\right)^{1+d/2} t^{-d/2}
\nn \\
& \le& c_4(T/r_0^2)^{1+d/2}\, t^{-d/2}.\label{e:nnw1122}
\end{eqnarray}
Thus by the monotonicity of the transition density and \eqref{e:nnw1122n} and \eqref{e:nnw1122} we obtain
\begin{eqnarray*}
 p_{D}(t, x, y)
&\le&c_{5} p(t, |x-y|/4)
\frac{\delta_D(x)}{\sqrt{t}}+ c_{6}  \left( t^{-d/2} \wedge  tj( |x-y|/4) \right)\frac{\delta_D(x)}{\sqrt{t}}\nn\\
&\le& \, c_{7} p(t, |x-y|/4)
\frac{\delta_D(x)}{\sqrt{t}} \\
&\leq &
c_7 p(t, |x-y|/4)
\left(\frac{\delta_D(x)}{\sqrt{t}} \wedge 1\right),
\end{eqnarray*}
where in second inequality \eqref{eq:globallbn} is used.

Combining these two cases establishes the claim \eqref{e:3.30}.
Thus by the semigroup property and the symmetry of $p_D(t, x, y)$ in $x$ and $y$, we
conclude from \eqref{e:3.30} that for every $t\in (0, T]$ and $x, y\in
D$,
\begin{eqnarray}
p_{D}(t,x,y)
&=& \int_{D} p_{D}({t}/2,x,z)p_{D}({t}/2,z,y)dz \nn\\
&\le&
c_7^2\left(1\wedge \frac{\delta_D(x)}{\sqrt{t}}\right)
\left(1\wedge \frac{\delta_D(
y)}{\sqrt{t}}\right) \int_{D}p(t/2, |x-z|/4) p({t}/{2}, |z-y|/{4})dz\nn\\
&\le&
c_7^2\left(1\wedge \frac{\delta_D(x)}{\sqrt{t}}\right)
\left(1\wedge \frac{\delta_D(
y)}{\sqrt{t}}\right) \int_{\R^d} p(t/2, |x-z|/4) p({t}/{2}, |z-y|/{4})dz\nn\\
&\le&
c_8 \left(1\wedge \frac{\delta_D(x)}{\sqrt{t}}\right)
\left(1\wedge \frac{\delta_D(
y)}{\sqrt{t}}\right)p(t, |x-y|/4). \nn
\end{eqnarray}

(ii) The lower bound in \eqref{e:nlu} is Theorem \ref{t:main}(iii). The proof of the
upper bound  in \eqref{e:nlu}
 is the same as that of Theorem \ref{t:main}(iv),
the only difference
is that we use part (i) of this theorem and  \eqref{globalupn}, instead of Theorem \ref{t:main}(ii), so that
$p_D(1, x, z) \leq c_1
(1\wedge {\delta_D(x)})(1\wedge {\delta_D(z)})
$ and $p_D(1, w, y) \leq c_1
(1\wedge {\delta_D(w)})(1\wedge {\delta_D(y)})
$
.
\qed

\section{Green function estimates}

In this section we give the proof of Corollary \ref{C:1.5}.

\medskip

\noindent{\bf Proof of Corollary \ref{C:1.5}:}
Put $T := \mbox{\rm diam}(D)^2$.
Recall that $g_D (x, y)$ is defined in \eqref{e:gD}.
By an argument similar to that for \cite[Corollary 1.2]{CKS},
one gets that (see the proofs of \cite[Theorem 5.0.8]{P} (for $d=1, 2$) and \cite[Theorem 6.2]{KS} (for $d \ge 3$) for details)
$$
 \int_0^T \left(1\wedge \frac{\delta_D(x)}{\sqrt{t}}\right)
\left(1\wedge \frac{\delta_D(
y)}{\sqrt{t}}\right)
 t^{-d/2}e^{-c_1|x-y|^2/t}  dt
+\int_T^\infty
  e^{-\lambda_1 t}\, \delta_D (x)\delta_D (y)  dt
\, \asymp\, g_D(x,y).
$$
Thus,
since $D$ is bounded, by Theorem \ref{t:main} (i) and (iii), we have
$ G_D(x,y)  \ge c_2  g_D(x,y)$,
which proves Corollary \ref{C:1.5}(i).

When the bounded $C^{1,1}$ open set $D$ satisfies \eqref{globalup}, by \eqref{e:jnearorigin}
and Theorem \ref{t:main} (ii) and (iv), we have
$$ G_D(x,y) \le  c_3 \left(g_D(x,y) +  \int_0^T\left( 1\wedge
\frac{\delta_D(x)}{\sqrt{t}}\right) \left( 1\wedge
\frac{\delta_D(y)}{\sqrt{t}}\right) \left( t^{-d/2}
\wedge \frac{t}{|x-y|^{d+2}}\right) dt \right).
$$
Therefore to prove Corollary \ref{C:1.5}(ii) it suffices to show that
\begin{eqnarray}
\int_0^T\left( 1\wedge
\frac{\delta_D(x)}{\sqrt{t}}\right) \left( 1\wedge
\frac{\delta_D(y)}{\sqrt{t}}\right) \left( t^{-d/2}
\wedge \frac{t}{|x-y|^{d+2}}\right) dt
\le c_4 g_D(x,y).\label{e:ch51}
\end{eqnarray}

 By the change of variable $u= \frac{|x-y|^2}{t}$, we have
\begin{eqnarray}
&& \int_0^T\left( 1\wedge
\frac{\delta_D(x)}{\sqrt{t}}\right) \left( 1\wedge
\frac{\delta_D(y)}{\sqrt{t}}\right) \left( t^{-d/2}
\wedge \frac{t}{|x-y|^{d+2}}\right) dt \nonumber\\
&=&\frac1{|x-y|^{d-2}}
 \int_{ |x-y|^2/T}^\infty
\left(u^{d/2-2} \wedge u^{-3} \right)  \left(1\wedge
\frac{ {\sqrt u} \delta_D(x)}{ |x-y| }\right)
\left(1\wedge \frac{ {\sqrt u} \delta_D(y)}{
|x-y| }\right) du.          \label{ID:5.1}
\end{eqnarray}
Since for every $x, y \in D$ and $r>0$,
\begin{equation}\label{e:ch5}
\left(1\wedge \frac{r \delta_D(x)} {|x-y|}\right)\, \left( 1\wedge
\frac{r \delta_D(y)} {|x-y|} \right) \, \le \, 1\wedge
\frac{r^2\delta_D(x)\delta_D(y)} {|x-y|^2},
\end{equation}
we have
\begin{eqnarray}
&&\frac1{|x-y|^{d-2}}  \int_{1}^\infty
\left(u^{d/2-2} \wedge u^{-3} \right) \left(1\wedge
\frac{ {\sqrt u} \delta_D(x)}{ |x-y| }\right)
\left(1\wedge \frac{ {\sqrt u} \delta_D(y)}{
|x-y| }\right) du   \nonumber\\
&=&
 \frac1{|x-y|^{d-2}}  \int_{1}^\infty
  u^{-2} \, \left(u^{-1/2} \wedge
\frac{   \delta_D(x)}{ |x-y| }\right)
\left(u^{-1/2} \wedge \frac{  \delta_D(y)}{
|x-y| }\right) du \nonumber\\
&\leq&
 \frac1{|x-y|^{d-2}}
\int_{1}^\infty  u^{-2} \,
  \left(1 \wedge
\frac{   \delta_D(x)}{ |x-y| }\right) \left(1
\wedge \frac{  \delta_D(y)}{
|x-y| }\right) du \nonumber\\
&\le & \frac1{ |x-y|^{d-2}}    \left(1\wedge \frac{
\delta_D(x)\delta_D(y)}{ |x-y|^2 }\right).
\label{e:ch2}
\end{eqnarray}

\medskip \noindent
 (1) When $d \ge 3$,  by \eqref{e:ch5}, we have
 \begin{eqnarray}
&&  \frac{1}{|x-y|^{d-2}} \int_{|x-y|^2/T}^{1 }
\left(u^{d/2-2} \wedge u^{-3} \right) \left( 1\wedge
\frac{ {\sqrt u} \delta_D(x)}{ |x-y| }\right)
\left(1\wedge \frac{ {\sqrt u} \delta_D(y)}{
|x-y| }\right)  du\nonumber\\
&\leq & \frac{1} {|x-y|^{d-2}}
 \left( 1\wedge \frac{  \delta_D(x)}{ |x-y|   }\right)
\left(1\wedge \frac{ \delta_D(y)}{ |x-y|   }
\right) \int_{0}^1  u^{d/2-2} du \nonumber \\
 &\leq &  \frac{2}{d-2} \frac{1} {|x-y|^{d-2}}
\left(1\wedge \frac{
\delta_D(x)\delta_D(y)}{ |x-y|^2 }\right) .  \label{e:ch4}
\end{eqnarray}
Combining \eqref{ID:5.1}, \eqref{e:ch2} and \eqref{e:ch4}, we arrive at \eqref{e:ch51} for $d \ge3$.
\medskip

For the other cases, we define
 \bee\label{e:ch6} u_0:=  \frac{
\delta_D(x) \delta_D(y)}{ |x-y|^{2} }.
 \eee
Clearly $1/u_0 \geq |x-y|^2/{\rm
diam}(D)^{2}=|x-y|^2/T$.

\medskip \noindent
(2) Suppose $d=2$. By \eqref{e:ch5} we have
 \begin{eqnarray}
&&  \frac{1}{|x-y|^{d-2}} \int_{|x-y|^2/T}^{1 }
\left(u^{d/2-2} \wedge u^{-3} \right) \left( 1\wedge
\frac{ {\sqrt u}\,  \delta_D(x)}{ |x-y|
}\right) \left(1\wedge \frac{ {\sqrt u}\, \delta_D(y)}{
|x-y| }\right)  du\nonumber\\
&\le&  \int_{|x-y|^2/T}^{1 } u^{-1}  \, \left(
1\wedge \frac{ u \,  \delta_D(x) \delta_D(y)}{
|x-y|^{2} }\right)   du \nonumber\\
&=& \int_{|x-y|^2/T}^{1 } u^{-1} {\bf 1}_{\{u\geq
1/u_0\}} du +  \int_{|x-y|^2/T}^{1 } u_0 {\bf 1}_{\{u<
1/u_0\}} du \nonumber \\
&=& \log (u_0 \vee 1) + u_0 \left( \frac1{u_0}\wedge 1 -
\frac{|x-y|^2}{T} \right) . \label{e:ch7}
\end{eqnarray}
Thus by  \eqref{ID:5.1}, \eqref{e:ch2} and \eqref{e:ch7},
\begin{eqnarray*}
&&\int_0^T\left( 1\wedge
\frac{\delta_D(x)}{\sqrt{t}}\right) \left( 1\wedge
\frac{\delta_D(y)}{\sqrt{t}}\right) \left( t^{-1}
\wedge \frac{t}{|x-y|^{4}}\right) dt \\
&\le & \left(1\wedge \frac{  \delta_D(x)}{
|x-y|   }\right) \left( 1\wedge \frac{
\delta_D(y)}{ |x-y|} \right) + \log (u_0 \vee
1) + u_0 \left( \frac1{u_0}\wedge 1 -
\frac{|x-y|^2}{T} \right)  \\
&\asymp& 1\wedge u_0 + \log (u_0 \vee 1) + u_0 \left( \frac1{u_0}\wedge
1 - \frac{|x-y|^2}{T} \right)   \\
&\asymp& 1\wedge u_0 + \log (u_0 \vee 1)
\asymp \log (1+ u_0 )  =  \log \left( 1+ \frac{  \delta_D(x) \delta_D
(y)}{ |x-y|^{2} }\right) .
\end{eqnarray*}
This proves \eqref{e:ch51} for $d =2$.

\medskip\noindent
(iii) Lastly we consider the case $d=1$. By
\eqref{e:ch5}  and \eqref{e:ch6},
\begin{eqnarray*}
&&  \frac{1}{|x-y|^{d-2}} \int_{|x-y|^2/T}^{1 }
\left(u^{d/2-2} \wedge u^{-3} \right) \left( 1\wedge
\frac{ {\sqrt u}\,  \delta_D(x)}{ |x-y|
}\right) \left(1\wedge \frac{ {\sqrt u}\, \delta_D(y)}{
|x-y| }\right)  du
\\
&\le& |x-y|
 \int_{|x-y|^2/T}^{1 } u^{-3/2}  \,
\left( 1\wedge \frac{ u \,  \delta_D(x)
\delta_D(y)}{
|x-y|^{2} }\right)   du
\\
&=& |x-y| \left(
\int_{|x-y|^2/T}^{1 } u^{-3/2}  {\bf
1}_{\{u\geq 1/u_0\}} du +  \int_{|x-y|^2/T}^{1 } u_0
u^{-1/2} {\bf 1}_{\{u<
1/u_0\}} du \right)
\\
 &=&  |x-y|\left( 2
 \left( (u_0\vee 1)^{1/2}-1\right) + 2 u_0 \left(
 (u_0\vee 1)^{-1/2} -
 \left(\frac{|x-y|^2}T \right)^{1/2}\right)
 \right).
\end{eqnarray*}
Thus by \eqref{ID:5.1}, \eqref{e:ch2} and
the last display, we have
\begin{eqnarray*}
&&\int_0^T\left( 1\wedge
\frac{\delta_D(x)}{\sqrt{t}}\right) \left( 1\wedge
\frac{\delta_D(y)}{\sqrt{t}}\right) \left( t^{-1/2}
\wedge \frac{t}{|x-y|^{3}}\right) dt  \\
&\le &
|x-y| \left(1\wedge u_0 \right)    + |x-y|\left( \left( (u_0\vee
1)^{1/2}-1\right) +   u_0 \left( (u_0\vee 1)^{-1/2}
-\left(\frac{|x-y|^2}T\right)^{1/2}\right)
\right) \\
&\asymp&
|x-y|\left(  u_0^{1/2} \wedge u_0 \right)
=  \left( \delta_D(x) \delta_D (y) \right)^{1/2}
\wedge  \frac{  \delta_D(x) \delta_D (y)}{ |x-y|} .
\end{eqnarray*}
This proves \eqref{e:ch51} for $d=1$. \qed

 \vskip 0.3truein

{\bf Zhen-Qing Chen}

Department of Mathematics, University of Washington, Seattle,
WA 98195, USA

E-mail: \texttt{zqchen@uw.edu}

\bigskip

{\bf Panki Kim}

Department of Mathematical Sciences and Research Institute of Mathematics,

Seoul National University, Building 27, 1 Gwanak-ro, Gwanak-gu Seoul 151-747,

Republic of Korea

E-mail: \texttt{pkim@snu.ac.kr}

\bigskip

{\bf Renming Song}

Department of Mathematics, University of Illinois, Urbana, IL 61801, USA

E-mail: \texttt{rsong@math.uiuc.edu}

\end{document}